\documentclass[11pt,openany,oneside]{amsart}
\usepackage[T1]{fontenc}
\usepackage[a4paper,left=2.5cm,right=2.5cm, top=3cm, bottom=3cm]{geometry}
\usepackage{amssymb}
\usepackage{enumitem}
\usepackage{lmodern}
\usepackage{hyperref}
\usepackage[hyperpageref]{backref} 
\usepackage{todonotes}

\numberwithin{equation}{section}
\newtheorem{Th}{Theorem}[section]
\newtheorem{Prop}[Th]{Proposition}
\newtheorem{Lem}[Th]{Lemma}
\newtheorem{Cor}[Th]{Corollary}

\newcounter{remarkcounter}

\newtheoremstyle{boldremark}   
{}                           
{}                           
{}                           
{}                           
{\bfseries}                 
{.}                          
{ }                          
{\thmname{#1}~\thmnumber{#2}} 

\theoremstyle{boldremark}
\newtheorem{RemTmp}{Remark}  

\newenvironment{Rem}
{\refstepcounter{remarkcounter} 
	\begin{RemTmp}} 
	{\end{RemTmp}}   

\newcommand{\N}{\mathbb{N}}
\newcommand{\R}{\mathbb{R}}
\newcommand{\C}{\mathbb{C}}

\newcommand{\cH}{{\mathcal H}}

\newcommand{\cP}{\mathcal{P}}
\newcommand{\cS}{{\mathcal S}}

\newcommand{\les}{\lesssim}
\newcommand{\divv}{{\rm div}}

 \def\dd{\, {\rm d}}

\newcommand{\tw}{\widetilde w}

\newcommand{\tJ}{\widetilde J}
\newcommand{\tK}{\widetilde K}
\newcommand{\tF}{\widetilde F}
\newcommand{\tc}{\widetilde c}
\newcommand{\tGamma}{\widetilde\Gamma}
\newcommand{\tgamma}{\widetilde\gamma}

\newcommand{\intd}{\int_{\R^2}}
\def\e{{\rm e}}
\def\ri{{\rm i}}

\newcommand{\Sra}{\cS_{r,a}}

\title[Positive solutions with prescribed mass for a planar Choquard equation with critical growth]{Positive solutions with prescribed mass for a planar Choquard equation with critical growth}
\date{\today}

\author[]{Ling Huang}%
\address[Ling Huang]{School of Mathematical Sciences, South China Normal University Guangzhou 510631, PR China}
\email{996987952@qq.com}

\author[]{Giulio Romani}%
\address[Giulio Romani]{Dipartimento di Scienza e Alta Tecnologia, Universit\`{a} degli Studi dell'Insubria and RISM-Riemann International School of Mathematics, Villa Toeplitz, Via G.B. Vico, 46 - 21100 Varese, Italy}
\email[Corresponding author]{giulio.romani@uninsubria.it}

\keywords{Normalised solution, Choquard equation, variational methods, critical exponential nonlinearities.}
\subjclass[2020]{35J20, 35J91, 35Q55, 35R09, 35B33.}

\begin{document}

\begin{abstract}
    We study normalised solutions for a Choquard equation in the plane with polynomial Riesz kernel and exponential nonlinearities, which are critical in the sense of Trudinger-Moser. For all prescribed values of the mass, we prove existence of a positive radial solution by a variational argument, which exploits a delicate analysis on the mountain pass level. Under an additional monotonicity assumption on the nonlinearity, such a solution turns out to be also a ground state in $H^1(\R^2)$. Our work extends the results by Dou, Huang, and Zhong in \cite{Dou2023} to the Choquard setting, improving in several directions those by Deng and Yu in \cite{DY}.
\end{abstract}
\maketitle

\section{Introduction}
We are interested in the existence of normalised solutions for the planar Choquard equation
\begin{equation}\label{eq}
    \begin{cases}
        -\Delta u +\lambda u = \left(I_\alpha\ast F(u)\right)f(u) &\text{ in } \R^2,\\
        \intd|u|^2\dd x=a^2.
    \end{cases}
\end{equation}
Here $I_\alpha(x):=\frac{c_{2,\alpha}}{|x|^{2-\alpha}}$ with $\alpha\in(0,2)$ denotes the (polynomial) Riesz kernel, and we assume an exponential critical growth for the positive and continuous nonlinearity $f$. Henceforth $F(t):=\int_0^tf(s)\dd s$. Since the value of the mass $a>0$ is prescribed, $\lambda\in\R$ appears as a Lagrange multiplier in the equation. Our main aim is to deal with problem \eqref{eq} in the limiting setting of the Sobolev embedding.
\vskip0.2truecm
For $N\geq2$ and $\alpha\in(0,N)$, denote the Riesz kernel by $I_{N,\alpha}(x):=\frac{c_{N,\alpha}}{|x|^{N-\alpha}}$, where $c_{N,\alpha}=\frac{\Gamma\left(\frac{N-\alpha}2\right)}{2^\alpha\pi^{\frac N2}\Gamma\left(\frac\alpha2\right)}$, and $\Gamma$ is the Gamma function. Choquard equations, namely Schr\"odinger equations with a nonlocal right-hand side, of the kind
\begin{equation}\label{Choq}
    -\Delta u +V(x)u = \left(I_{N,\alpha}\ast F(u)\right)f(u)\quad\text{ in } \R^N,
\end{equation}
where $V$ is an external potential, are equivalent formulations of the Schr\"odinger-Poisson systems
\begin{equation}\label{SP:system}
    \begin{cases}
        -\Delta u +V(x)u = \Phi f(u) &\text{ in } \R^N,\\
        (-\Delta)^{\frac\alpha2}\Phi = F(u) &\text{ in } \R^N,
    \end{cases}
\end{equation}
since the kernel $I_{N,\alpha}$ is the fundamental solution of the operator $(-\Delta)^{\frac\alpha2}$ in $\R^N$. For this reason, they emerge in various fields of physics, e.g. in Hartree theory for crystals, astrophysics, electromagnetism, and quantum mechanics. We refer the interested reader to \cite{BF,LRZ} for the background. There is a vast literature on Choquard equations with polynomial kernels, which flourished especially in the last two decades, since in this context new intriguing phenomena appear: we refer to the important works \cite{MVS,MVS2,CZ} when $N\geq3$, and to the recent developments \cite{BvS,ACTY,AFS,R2} in the planar case $N=2$. Note that, as in the Schr\"odinger case \cite{dFMR}, an exponential growth of the nonlinearity $f$ in \eqref{eq} is allowed thanks to Ruf's inequality, which extends to the whole $\R^2$ the classical results by Trudinger and Moser for bounded domains. This introduces additional difficulties related to the lack of compactness.

Equations of the kind \eqref{SP:system} derive from time-dependent Schr\"odinger-Poisson systems
\begin{equation}\label{SP:system_t}
    \begin{cases}
        \ri\frac{\dd\Psi}{\dd t}=\Delta\Psi + \omega f(\Psi) &\text{ in } \R^N,\\
        (-\Delta)^{\frac\alpha2}\omega = F(\Psi) &\text{ in } \R^N,
    \end{cases}
\end{equation}
where $\Psi:\R\times\R^N\to\C$ is the time-dependent wave function and $\omega:\R\times\R^N\to\C$ has the role of an internal potential which takes into account nonlocal self-interactions of the wave function. Indeed, under the time-harmonic ansatz $(\Psi(t,x),\,\omega(t,x))=\e^{-\ri\lambda t}(u(x),\,\Phi(x))$ with $\lambda\in\R$, if $f(\e^{-\ri\lambda t}s)=\e^{-\ri\lambda t}f(s)$ for all $\lambda,s\in\R$ and $t>0$, and similarly for $F$, then \eqref{SP:system_t} reduces to \eqref{SP:system} with $V(x)\equiv\lambda$. Since it is well-known that in this case $\|\Psi(\cdot,t)\|_2=\|\Psi(\cdot,0)\|_2=\|u\|_2$ for all $t>0$, then it makes sense to look for solutions of the stationary Choquard equation with a prescribed mass. This, in dimension two, leads to solve \eqref{eq}, where $\lambda$ is unknown and appears as a Lagrange multiplier. This constrained problem is also interesting from a mathematical point of view, since the usual methods to investigate the unconstrained equation \eqref{Choq} have to be adapted, and in particular new behaviours and phenomena appear. In fact, the solutions of \eqref{eq}, called \textit{normalised}, correspond to the critical points of the functional 
\begin{equation}\label{J}
    J(u)=\frac12\intd|\nabla u|^2\dd x-\frac12\intd\left(I_\alpha\ast F(u)\right)F(u)\dd x
\end{equation}
constrained to the $L^2$-sphere
$$S_a:=\{u\in H^1(\R^2)\,|\,\|u\|_2=a\}\,,$$
and the nature of critical points depends on geometry of $J$ on $S_a$. In the Schr\"odinger case in $\R^N$, that is, when the nonlocal right-hand side in \eqref{eq} is replaced by a local nonlinearity $f(u)$, the corresponding functional is bounded from below provided the growth of $F$ does not exceed $t^{2+\frac4N}$, and one can look for minima for $J$ on $S_a$, see \cite{Shibata}. This case is known as \textit{mass-subcritical}. On the other hand, namely in the \textit{mass-supercritical} case, when the growth of $F$ is still Sobolev subcritical but higher than the above threshold, then $\displaystyle{\inf_{S_a}}\,J=-\infty$, and one has to look for critical points with a minimax characterisation. We remark that in the latter case the boundedness of Palais-Smale sequences is not assured. In \cite{Jeanjean1997}, Jeanjean managed to overcome this problem by applying Ekeland's variational principle to an auxiliary functional obtained as the composition of $J$ with a fiber map which preserves the $L^2$-norm: this permits to obtain a (PS)-sequence with the further properties that the Poho\v zaev functional is almost satisfied, which is a striking property to show the boundedness of the sequence. This strategy has become nowadays well-known and several refinements and extension of the results by Jeanjean may be found, see e.g. \cite{HT,JL,BM,JZZ}. For the corresponding Choquard problem \eqref{eq} in $\R^N$, the presence of the nonlocal term does not affect the occurrence of the phenomena just described, except for a different threshold defining the mass-critical case, which in this case appears when $F$ behaves like $t^{\frac{N+\alpha+2}N}$, see e.g. \cite{BJL,LY,YCT,BLL,XZ}.
\vskip0.2truecm
More recently, the above mentioned results for the Schr\"odinger case have been extended to the planar setting, where the nonlinearity may be supposed to be exponentially growing, see \cite{AJM,CLY,Dou2023}. The case of nonlinearities with critical exponential growth in the sense of Trudinger-Moser, namely behaving at infinity like $\e^{\gamma t^2}$ for some $\gamma>0$, corresponds to the critical Sobolev case, since phenomena of lack of compactness arise. Therefore, one is compelled to recover compactness by means of fine estimates on suitable energy levels. It is common in the literature, and so also in \cite{AJM,CLY}, to prescribe a strong growth condition of the kind
\begin{equation}\label{global_bound_below_Schr}
    F(t)\geq\nu|t|^\mu\qquad\mbox{for all}\ \,t>0\,,
\end{equation}
with $\mu>\mu_0$, for a suitable $\mu_0$, and a constant $\nu$ so large to force a compression of the mountain pass level under the first noncompactness level of the functional; usually a bound from below of $\nu$ is given but it is often not explicit. Despite the fact that \eqref{global_bound_below_Schr} is clearly satisfied for large values of $t$, since $F$ is assumed to have an exponential growth at infinity, however, it is not of practical verification in concrete cases in a neighbourhood of $0$, due to the large value of $\nu$. In \cite{CRTY,Dou2023} this assumption has been removed via a delicate analysis on Moser-type sequences and substituted by a growth condition at infinity in the spirit of \cite{dFMR}, which is satisfied in most of the cases.
\vskip0.2truecm
To our knowledge, the only work about normalised solutions for the planar Choquard equation of the form \eqref{eq} is \cite{DY}, where both cases of subcritical and critical exponential nonlinearities in the sense of Trudinger-Moser are considered. If, on the one hand, the analysis in the subcritical case is complete, for the critical case the authors rely on the assumption from below \eqref{global_bound_below_Schr} and, moreover, existence of a normalised solution is proved only for small masses $a\in(0,a_0)$, where $a_0$ is explicit and depends on the growth of $f$. We also refer to \cite{ChenWang, ChenWang2} for the higher-order case in the conformal dimension. 
\vskip0.2truecm
In this manuscript, we aim, on the one hand, to extend the results in \cite{Dou2023} to the nonlocal Choquard case, and, on the other hand, to improve those in \cite{DY} in several directions: first, by means of a condition at infinity in the spirit of \cite{dFMR}, we remove the assumption from below \eqref{global_bound_below_Schr}; moreover, we manage to prove existence of a positive radial normalised solution for \eqref{eq} for \textit{all} values of the mass in the critical case which, under a monotonicity assumption on $f$, also turns out to be a ground state solution; finally, we provide an estimate for the range of the frequency $\lambda$ which depends on all the structural constants into play. Incidentally, we also weaken the growth condition of the nonlinearity at $0$ appearing in \cite{DY}, which is not optimal for the Choquard case.
\vskip0.2truecm

From now on, we consider nonlinearities $f$ satisfying the following set of assumptions.
\begin{enumerate}
    \item[($f_0$)] $f\in C(\R,\R^+)$, $f(t)=0$ for all $t\leq0\,$;
    \item[($f_1$)] $f(t)=o(t^{1+\frac\alpha2})$ as $t\to0^+$;
    \item[($f_2$)] there exists $\gamma_0>0$ such that $$\lim_{t\to+\infty}\frac{f(t)}{\e^{\gamma t^2}}=\begin{cases}
		  0 &\text{ if } \gamma>\gamma_0\,,\\
		+\infty &\text{ if } \gamma<\gamma_0\,;
	\end{cases}$$
    \item[($f_3$)] there exists $\mu>2+\frac\alpha2$ such that $0<\mu F(t)\leq f(t)t$ for all $t>0\,$;
    \item[($f_4$)] $\displaystyle{\lim_{t\to+\infty}}\frac{F(t)}{f(t)t}=0\,$;
    \item[($f_5$)] there exists $\beta_0>0$ such that $\displaystyle{\lim_{t\to+\infty}}\frac{tf(t)}{\e^{\gamma_0t^2}}\geq\beta_0\,$.
    
\end{enumerate}

\begin{Rem}
    We use ($f_0$) in order to retrieve positivity for the solutions; ($f_1$) and ($f_2$) are the growth conditions at $0$ and $+\infty$, respectively, which sort the problem as \textit{mass-supercritical}; ($f_3$) is the usual Ambrosetti-Rabinowitz condition, and we note that the nonlocal structure of the equations enters here and in ($f_1$); $(f_4)$ plays a role for compactness; ($f_5$) is used to properly estimate the mountain pass level and is the same condition as in \cite{dFMR}, but here there is no need for $\beta_0$ to be a large parameter. A more detailed description of the assumptions, as well examples of nonlinearities which fulfil them, will be given in Section \ref{Sec_Prel}.
\end{Rem}

We are now in a position to present our main results. In the following, $H^1_r(\R^2)$ denotes the subspace of $H^1(\R^2)$ composed by radial functions.
\begin{Th}\label{Thm}
    Suppose conditions ($f_0$)-($f_5$) are satisfied and define $\Sra:=\cS_a\cap H_r^1(\R^2)$. Then, for all $a>0$, problem \eqref{eq} has a mountain pass type solution $(u_a,\lambda_a)\in\Sra\times\R^+$. Moreover, $\lambda_a\in\left(0,\frac{(2+\alpha)^2\pi}{2\gamma_0a^2\left(\mu-\left(2+\frac\alpha2\right)\right)}\right)$ and $u_a>0$.
\end{Th}
\begin{Rem}
    Note that the upper bound on $\lambda_a$ depends on all the structural constant of the problem: the critical exponential growth by $\gamma_0$, the mass by $a^2$, the Riesz kernel by $\alpha$, and the Ambrosetti-Rabinowitz constant by $\mu>2+\frac\alpha2$.
\end{Rem}
Under an additional monotonicity assumption, we are also able to show that the mountain pass type solution found in Theorem \ref{Thm} is indeed a normalised ground state solution. We recall that a normalised solution $u$ is a \textit{ground state} if $J(u)$ is minimal among all normalised solutions.
\begin{Th}\label{Thm_gs}
    In addition to the assumptions of Theorem \ref{Thm}, suppose that
    \begin{itemize}
        \item[($f_6$)] the function $\tF(t):=f(t)t-\frac{2+\alpha}2F(t)$ is nondecreasing in $(0,+\infty)\,.$
    \end{itemize}
    Then the solution $u_a$ obtained in Theorem \ref{Thm} is a normalised ground state solution of \eqref{eq}.
\end{Th}
\begin{Rem}
    We stress the fact that, although we work in $H^1_r(\R^2)$ to prove the existence of the solution $u_a$ in Theorem \ref{Thm}, it turns out that $u_a$ is of minimal energy in the whole $H^1(\R^2)$ and not only in its radial subspace $H^1_r(\R^2)$. This is accomplished by an argument which uses rearrangement inequalities, see Lemma \ref{c2=c2rad} in Section \ref{sec_proofs}.
\end{Rem}

\vskip0.2truecm
\begin{Rem}
    We are confident that, by a direct extension of our arguments, and under a suitable but straightforward modification of the assumptions ($f_1$)-($f_6$), one may also solve the corresponding $N$-dimensional problem in the limiting Sobolev setting, namely
    \begin{equation*}
        \begin{cases}
            -\Delta_N u +\lambda|u|^{N-2}u = \left(I_{\alpha,N}\ast F(u)\right)f(u) &\text{ in } \R^N,\\
            \int_{\R^N}|u|^N\dd x=a^N,
        \end{cases}
    \end{equation*}
    where $\Delta_N u:=\divv(|\nabla u|^{N-2}\nabla u)$ is the $N$-Laplace operator, $\alpha\in(0,N)$, and $a>0$. Here, for the sake of a more transparent exposition, we just deal with the planar case.
\end{Rem}

\begin{Rem}
    We also point out that, in dimension two, in order to keep the connection with \textit{local} Schr\"odinger-Poisson systems, namely \eqref{SP:system} with $\alpha=2$, the kernel $I_\alpha$ should be replaced by $I_2=\tfrac1{2\pi}\log\frac1{|\cdot|}$. However, its unboundedness from above and below adds further intricacies, in particular at the level of detecting an appropriate functional settings, see \cite{CT,LRTZ,CDL,R1,BRT} for Choquard equations in this "double limiting" setting. We plan to address this problem in a future work.
\end{Rem}
\vskip0.2truecm
Before going into the details of the proofs, let us briefly describe our method. The problem is variational, in the sense that normalised solutions of \eqref{eq} may be found by looking at critical points of the functional $J$ defined in \eqref{J} constrained to the $L^2$-sphere $\mathcal{S}_a$. Since the problem is autonomous, in order to retrieve compactness we may equivalently work in $H^1_r(\R^2)$, since here we can benefit from the compactness of the embedding $H^1_r(\R^2)\hookrightarrow L^q(\R^2)$ for all $q>2$, see \cite{Lions82}. Indeed, by Palais' principle of symmetric criticality, see \cite{Palais1979}, it is standard to prove that the solutions in $H_r^1(\R^2)$ are in fact solutions in the whole $H^1(\R^2)$. Thanks to ($f_1$)-($f_2$), the functional $J:H^1_r(\R^2)\to\R$ is then well-defined and $C^1$ by means of Cao's and Hardy-Littlewood-Sobolev's inequalities (Propositions \ref{Cao_ineq} and \ref{HLS} below, respectively). Moreover, the mass-supercritical nature of the problem yields a mountain pass type geometry (Lemma \ref{MP_geometry}). After \cite{Jeanjean1997,BLL}, it is well-known that this implies the existence of a Palais-Smale-Poho\v zaev sequence ((PSP)-sequence for short) at the mountain pass level $c_{mp}$, namely a Palais-Smale sequence $\{u_n\}$ at level $c_{mp}$ such that $\cP(u_n)\to0$ as $n\to+\infty$, where $\cP$ is the Poho\v zaev functional associated to \eqref{eq}, see \eqref{Pohozaev}. The latter property is striking in order to prove the boundedness of the sequence $\{u_n\}$. A delicate part of the analysis is to get the desired upper bound on the mountain pass level $c_{mp}$. This will be obtained in Lemma \ref{20240605-L1} by exploiting assumption ($f_5$) and using Moser sequences. This estimate is crucial to prove the strong convergence of the (PSP)$_{c_{mp}}$-sequence to a nontrivial critical point of $J$ on $S_a$ in Section \ref{sec_PSP}. Finally, to prove that the mountain pass solution we find is also a ground state solution, we rely on the monotonicity assumption ($f_6$) and show that  the functional evaluated on some fiber paths reaches its maximum in a unique point which belongs to the Poho\v zaev manifold $\{\cP(u)=0\}$. Once this geometry is proved, it is then plain to compare the mountain pass level and the least-energy level.
\vskip0.2truecm
\paragraph{\textbf{Overview.}} After a brief Section \ref{Sec_Prel} in which we recall some useful results needed in particular to deal with the exponential growth of our nonlinearity, in Section \ref{Sec_MP} we check that the functional $J\big|_{S_a}$ (similarly for $J\big|_{\Sra}$) possesses a mountain pass structure, we prove the existence of a $(PSP)$-sequence at level $c_{mp}$, as well as a suitable upper bound of the mountain pass level $c_{mp}$. In Section \ref{sec_PSP}, we show that the $(PSP)_{c_{mp}}-$ sequence is bounded and has a weak limit, deferring the proof of Theorems \ref{Thm} and \ref{Thm_gs} to Section \ref{sec_proofs}.
\vskip0.2truecm
\paragraph{\textbf{Notation.}} For $R>0$ and $x_0\in\R^2$ we denote by $B_R(x_0)$ the ball of radius $R$ and centre $x_0$. Given a set $\Omega\subset\R^N$, we denote $\Omega^c:=\R^N\setminus\Omega$, its Lebesgue measure by $|\Omega|$, and its characteristic function by $\chi_\Omega$. For $p\in[1,+\infty]$ the Lebesgue space of $p$-integrable functions is denoted by $L^p(\R^N)$ with norm $\|\cdot\|_p$. For $\tau>1$ we define its conjugate H\"older exponent as $\tau':=\frac \tau{\tau-1}$. The symbol $\lesssim$ indicates that an inequality holds up to a multiplicative constant depending only on structural constants. Finally, $\text{o}_n(1)$ denotes a vanishing real sequence as $n\to+\infty$. Hereafter, the letter $C$ will be used to denote positive constants which are independent of relevant quantities and whose value may change from line to line.

\section{Preliminary results}\label{Sec_Prel}

\begin{Prop}[Cao's inequality, \cite{Cao}]\label{Cao_ineq}
    (i) If $\gamma>0$ and $u\in H^1(\R^2)$, then
    $$\intd\left(\e^{\gamma u^2}-1\right)\dd x<+\infty\,.$$
   (ii) Moreover, if $\|\nabla u\|_2\leq1$, $\|u\|_2\leq M$ and $\gamma<4\pi$, then there exists a constant $C=C(M,\gamma)$ such that
    $$\intd\left(\e^{\gamma u^2}-1\right)\dd x<C(M,\gamma)\,.$$
\end{Prop}

Besides the simple inequality
\begin{align*}
    (\e^s-1)^t\leq\e^{ts}-1,\quad\text{for}~t>1~\text{and}~s\geq0\,.
\end{align*}
which will be frequently used in the following, we shall recall a useful result which will be important in the convergence analysis in Section \ref{sec_PSP}.

\begin{Lem}\label{lemma 2.2}
    Let $u_0\in H^1(\R^2)$, and $\{u_n\}$ be a bounded sequence in $H^1(\R^2)$, such that
    $$\limsup_{n\to+\infty}\|\nabla(u_n-u_0)\|_2^2<\frac{(2+\alpha)\pi}{\gamma_0}\,.$$
    Then for $\gamma>\gamma_0$ close to $\gamma_0$ the sequence $\Big\{\e^{\frac{4\gamma}{2+\alpha} u_n^2}-1\Big\}$ is bounded in $L^t(\R^2)$ provided $t>1$ is close to $1$.
\end{Lem}
\begin{proof}
  We follow the line of \cite[Lemma 2.8]{Zhang2022}, see also \cite[Lemma 2.4]{Dou2023}. For $\gamma>\gamma_0$ close to $\gamma_0$, and $t>1$ close to 1, we still have that
  \begin{align}\label{2024/12/17/e1}
     \limsup_{n\to+\infty}\left(\frac{4\gamma t}{2+\alpha}\|\nabla(u_n-u_0)\|_2^2\right) <4\pi.
  \end{align}
 Moreover, by using the simple identity
 \begin{equation}\label{exp_identity}
    \e^{a+b}-1=(\e^a-1)(\e^b-1)+(\e^a-1)+(\e^b-1)\,,
 \end{equation}
 and choosing $\sigma>0$ small enough and $\eta>1$, by H\"older's inequality we infer
 \begin{equation}\label{lemma2.2_proof}
    \begin{split}
        \intd\left(\e^{\frac{4\gamma}{2+\alpha}u_n^2}-1\right)^t\dd x&\leq\intd\left(\e^{\frac{4\gamma t}{2+\alpha}u_n^2}-1\right)\dd x\\
         &=\intd\left(\e^{\frac{4\gamma t}{2+\alpha}(u_n-u_0+u_0)^2}-1\right)\dd x\\
         &\leq\intd\left(\e^{\frac{4\gamma t}{2+\alpha}\left[(1+\sigma)(u_n-u_0)^2+(1+\frac{1}{\sigma})u_0^2\right]}-1\right)\dd x\\
         &\leq\left(\intd\Big(\e^{\frac{4\gamma t(1+\frac1\sigma)\eta'}{2+\alpha}u_0^2}-1\Big)\dd x\right)^\frac1{\eta'}\!\left(\intd\Big(\e^{\frac{4\gamma t(1+\sigma)\eta}{2+\alpha}|u_n-u_0|^2}-1\Big)\dd x\right)^\frac1\eta\\
         &\quad+\intd\left(\e^{\frac{4\gamma t}{2+\alpha}(1+\sigma)(u_n-u_0)^2}-1\right)\dd x+\intd\left(\e^{\frac{4\gamma t}{2+\alpha}\left(1+\frac{1}{\sigma}\right)u_0^2}-1\right)\dd x\,.
    \end{split}
 \end{equation}
 By choosing $\eta>1$ close to 1, \eqref{2024/12/17/e1}, and the uniform estimate in Proposition \ref{Cao_ineq} (ii), then
 $$\intd\Big(\e^{\frac{4\gamma t(1+\sigma)\eta}{2+\alpha}|u_n-u_0|^2}-1\Big)\dd x\leq C$$
 for all $n$ large, while
 $$\intd\Big(\e^{\frac{4\gamma t(1+\frac1\sigma)\eta'}{2+\alpha}u_0^2}-1\Big)\dd x<+\infty$$
 by Proposition \ref{Cao_ineq} (i). The last two terms in \eqref{lemma2.2_proof} can be handled in the same way.
\end{proof}

Throughout the paper, we will make large use of the well-known Hardy-Littlewood-Sobolev inequality (see \cite[Theorem 4.3]{LiebLoss}).
\begin{Prop}[HLS inequality]\label{HLS}
    Let $N\geq1$, $\mu\in(0,N)$, and $s,r>1$ with $\tfrac1s+\tfrac\mu N+\tfrac1r=2$. There exists a constant $C=C(N,\mu,s,r)$ such that for all $f\in L^s(\R^N)$ and $h\in L^r(\R^N)$ one has
    $$\int_{\R^N}\left(\frac1{|\cdot|^\mu}\ast f\right)\!h\dd x\leq C\|f\|_s\|h\|_r\,.$$
\end{Prop}
\vskip0.2truecm
\paragraph{\textbf{Consequences of the assumptions.}}
To end this Section, let us point out some immediate consequences of ($f_0$)-($f_6$) which will be of use in our analysis, together with some comments in this regard.
\begin{enumerate}
    \item[i)] By $(f_1)$ and $(f_2)$, for any $q>1$, $\gamma>\gamma_0$, and $\varepsilon>0$, there exists $C_\varepsilon>0$ such that
    \begin{equation}\label{estf1}
        |f(t)|\leq\varepsilon|t|^{1+\frac\alpha2}+C_\varepsilon|t|^{q-1}\left(\e^{\gamma t^2}-1\right),\quad\forall~ t\in\R
    \end{equation}
    and
    \begin{equation}\label{estF2}
        |F(t)|\leq\varepsilon|t|^{2+\frac\alpha2}+C_\varepsilon|t|^q\left(\e^{\gamma t^2}-1\right), \quad\forall~ t\in\R\,.
    \end{equation}
    \item[ii)] $(f_4)$ is a condition at $\infty$ which is used to gain compactness. It is weaker than the assumption
    \begin{itemize}
        \item[($f_4'$)] There exist $s_0$, $L_0>0$ such that $F(s)\leq L_0 f(s)$ for $s \geq s_0\,$,
    \end{itemize}
    which is classical for problems dealing with the exponential critical growth. Moreover, from $(f_1)$ and $(f_4)$ we also infer that for any $\varepsilon>0$, we can find $C_\varepsilon>0$ such that
    \begin{equation}\label{f1f4}    
        F(t)\leq \varepsilon f(t)t+C_\varepsilon t^{2+\frac\alpha2}\qquad\mbox{for all}\ \,t\in\R\,.
    \end{equation}
    \item[iii)] ($f_5$) is the well-known de Figueiredo-Miyagaki-Ruf condition \cite{dFMR} and is crucial in order to estimate the mountain pass level, see Lemma \ref{20240605-L1} below. We note here that such an assumption avoids the prescription of a lower bound on $F$ like \eqref{global_bound_below_Schr} as in \cite{DY}, with $\mu>2+\frac\alpha 2$ and with $\nu$ \textit{large enough}: although widely used in the literature, the latter is not of practical verification. We remark that the choice of a large $\nu$ is essential in their paper, as it allows for a control on the mountain pass level.
    \item[iv)] From ($f_5$) one may easily infer that
    \begin{equation}\label{F_Ruf}
       \lim_{s\to+\infty}\frac{F(s)s^2}{\e^{\gamma_0s^2}}\geq\frac{\beta_0}{2\gamma_0}\,,
    \end{equation}
    see \cite[Lemma 3.4]{Dou2023}. Indeed, by de l'H\^opital's rule,
    \begin{equation*}
        \lim_{s\to+\infty}\frac{F(s)s^2}{\e^{\gamma_0s^2}}=\lim_{s\to+\infty}\frac{f(s)}{-2s^{-3}\e^{\gamma_0s^2}+2\gamma_0s^{-1}\e^{\gamma_0s^2}}=\frac1{2\gamma_0}\lim_{s\to+\infty}\frac{sf(s)}{\e^{\gamma_0s^2}}\geq\frac{\beta_0}{2\gamma_0}\,.
    \end{equation*}
    \item[v)] An example of a nonlinearity which satisfies assumptions ($f_1$)-($f_6$) is
    \begin{equation*}
        f(s)=s^{\sigma-1}\chi_{\{0<s<s_0\}}(s)+\frac{\beta_0(\gamma_0s^2-1)\e^{\gamma_0s^2}}{\gamma_0s^3}\chi_{\{s\geq s_0\}}(s)\,,
    \end{equation*}
    where $\sigma\in(2+\alpha,6)$ and $s_0$ is chosen so that $f$ is continuous on $\R^+$. Indeed, ($f_1$), ($f_2$), and ($f_5$) are of immediate verification; for $s>s_0$ we have
    $$F(s)=\frac{s^\sigma}\sigma+\frac{\beta_0}{2\gamma_0}\left(\frac{\e^{\gamma_0s^2}}{s^2}-\frac{\e^{\gamma_0s_0^2}}{s_0^2}\right),$$
    and
    $$\widetilde F(s)=2\beta_0\e^{\gamma_0s^2}\gamma_0s^3\left(\gamma_0^2s^4-\gamma_0\left(1+\frac{2+\alpha}4\right)s^2+\left(1+\frac{2+\alpha}4\right)\right)\,,$$
    from which ($f_4$) and ($f_6$) hold, respectively. Finally, if $\sigma>6$ it is easy to show that $f(s)/s^{\sigma-1}$ is nondecreasing on $\R^+$, which is equivalent to $f'(s)s\geq(\sigma-1)f(s)$. This implies that the function $g(s):=\sigma F(s)-sf(s)$ is nonpositive, which is ($f_3$) with $\mu=\sigma>2+\alpha$.
    
    The same can be proved also for nonlinearities $f$ for which $F$ is of the kind
    $$F(s)=s^p\chi_{\{0<s<s_0\}}(s)+B\frac{\e^{\gamma_0s^2}}{s^q}\chi_{\{s\geq s_0\}}(s)$$
    with $p>2+\frac\alpha2$, $\gamma_0>0$, $q\leq 2$ (in order for ($f_5$) to be satisfied), $s_0>\sqrt{\frac{\max\{p+q,0\}}{2\gamma_0}}$ so that both ($f_3$) and ($f_6$) are fulfilled, and $B\in\R^+$ so that $f\in C(\R)$.
\end{enumerate}

\section{The mountain pass structure}\label{Sec_MP}

We first verify that the functional $J$ is well-defined in $H^1(\R^2)$ and that on the sphere $\Sra$ possesses a mountain pass type geometry. To this aim, for $u\in H^1(\R^2)$, we introduce the $L^2$-invariant scaling
$$\cH(u,s):=\e^su(\e^s\cdot)\,,\qquad s\in\R\,,$$
and the functional $\tJ$ defined by
\begin{equation}\label{tJ}
    \tJ(u,s):=J(\cH(u,s))=\frac{\e^{2s}}2\intd|\nabla u|^2\dd x-\frac1{2\e^{(2+\alpha)s}}\intd\left(I_\alpha\ast F(\e^su)\right)F(\e^su)\dd x\,.
\end{equation}
We first observe that $\tfrac{\dd\tJ(u,s)}{\dd s}\big|_{s=0}=\cP(\cH(u,s))$, where
\begin{equation}\label{Pohozaev}
    \cP(u):=\intd|\nabla u|^2\dd x+\frac{2+\alpha}2\intd\left(I_\alpha\ast F(u)\right)F(u)\dd x-\intd\left(I_\alpha\ast F(u)\right)f(u)u\dd x
\end{equation}
is the Poho\v zaev functional associated to \eqref{eq}. We recall that if $(\lambda,u)$ is any couple weakly solving problem \eqref{eq}, then $u\in\cP(a)$, where
\begin{equation}\label{Pohozaev_manifold}
    \cP(a)=\{u\in\cS_a\,|\,\cP(u)=0\}\,.
\end{equation}
Indeed, testing the equation with the solution itself one gets
$$\|\nabla u\|_2^2+\lambda\|u\|_2^2=\intd\left(I_\alpha\ast F(u)\right)f(u)u\dd x$$
which, combined with the Poho\v zaev identity relative to \eqref{eq}, that is
\begin{align}\label{20240531-e8}
    \lambda\|u\|_2^2=\left(1+\frac\alpha2\right)\int_{\R^2}(I_\alpha\ast F(u))F(u)\dd x\,,
\end{align}
yields $\cP(u)=0$. We point indeed out, that the Poho\v zaev identity for \eqref{eq} set in $\R^N$ with $N\geq2$ is
	$$\frac{N-2}2\intd|\nabla u|^2\dd x+\frac{\lambda N}2\intd|u|^2\dd x=\frac{N+\alpha}2\intd\left(I_\alpha\ast F(u)\right)F(u)\dd x\,,$$
	which reduces to \eqref{20240531-e8} in the case $N=2$. Such an identity can be deduced by the same argument used in \cite[Proposition 3.1]{MVS2}, which deals with the simpler case $F(u)=|u|^p$.

\vskip0.2truecm
Although the next step is rather standard, we do it plainly for the sake of completeness, to show the role of assumptions $(f_1)-(f_2)$, as well as to address some inaccuracies in \cite{DY}.
\begin{Lem}\label{H_J}
    Under $(f_0)-(f_3)$, the functional $J$ is well-defined in $H^1(\R^2)$. Moreover, for a fixed $u\in\Sra$, one has
    \begin{enumerate}
        \item[(i)] $\|\nabla\cH(u,s)\|_2\to0$ and $\tJ(u,s)\to0$ as $s\to-\infty\,$;
        \item[(ii)] $\|\nabla\cH(u,s)\|_2\to+\infty$ and $\tJ(u,s)\to-\infty$ as $s\to+\infty\,$.
    \end{enumerate}
\end{Lem}
\begin{proof}
    We focus on the second term of $J$ and we apply the Hardy-Littlewood-Sobolev inequality (Proposition \ref{HLS}) with $N=2$ and $\mu=2-\alpha$ and $r=s=\tfrac4{2+\alpha}$, and \eqref{estF2}, obtaining
    \begin{equation}\label{estimate_convolution_term}
        \begin{split}
            \bigg|\intd&\left(I_\alpha\ast F(u)\right)F(u)\dd x\bigg|\les\|F(u)\|_{\frac 4{2+\alpha}}^2\\
            &\les\varepsilon\left(\intd|u|^{\frac{2(\alpha+4)}{\alpha+2}}\dd x\right)^\frac{2+\alpha}2\!+C_\varepsilon\left(\intd|u|^{\frac{4qt'}{2+\alpha}}\dd x\right)^\frac{2+\alpha}{2t'}\!\left(\intd\left(\e^{\frac{4\gamma t}{2+\alpha}u^2}-1\right)\dd x\right)^\frac{2+\alpha}{2t}\!.
        \end{split}
    \end{equation}
    Since $\frac{2(\alpha+4)}{(\alpha+2)}>2$, the first term is well-defined. The same holds for the second term, since one may choose $q>1$ arbitrary large (up to a larger constant $C_\varepsilon$). The third term is finite for all fixed $u\in H^1(\R^2)$ thanks to Proposition \ref{Cao_ineq} (i).

    Since $\|\nabla\cH(u,s)\|_2=\e^{s}\|\nabla u\|_2$, the properties (i)-(ii) for $\cH(u,s)$ follow immediately. Moreover, by \eqref{estimate_convolution_term}, one has
    \begin{equation}\label{estimate_convolution_term_H}
        \begin{split}
            &\bigg|\intd\left(I_\alpha\ast F(\cH(u,s))\right)F(\cH(u,s))\dd x\bigg|\\
            &\les\varepsilon\left(\intd|\cH(u,s)|^{\frac{2(\alpha+4)}{\alpha+2}}\dd x\right)^\frac{2+\alpha}2\!+C_\varepsilon\left(\intd|\cH(u,s)|^{\frac{4qt'}{2+\alpha}}\dd x\right)^\frac{2+\alpha}{2t'}\\
            &\quad\times\left(\intd\left(\e^{\frac{4\gamma t}{2+\alpha}\|\nabla\cH(u,s)\|_2^2\left(\frac{\cH(u,s)}{\|\nabla\cH(u,s)\|_2}\right)^2}-1\right)\dd x\right)^\frac{2+\alpha}{2t}\\
            &\les\varepsilon\left(\intd|\cH(u,s)|^{\frac{2(\alpha+4)}{\alpha+2}}\dd x\right)^\frac{2+\alpha}2\!\!+C_\varepsilon\left(\intd|\cH(u,s)|^{\frac{4qt'}{2+\alpha}}\dd x\right)^\frac{2+\alpha}{2t'},
        \end{split}
    \end{equation}
    where in the last inequality we used the uniform version of Cao's inequality given by Proposition \ref{Cao_ineq} (ii), provided one chooses $s$ sufficiently negative so that $\frac{4\gamma t}{2+\alpha}\|\nabla\cH(u,s)\|_2^2=\frac{4\gamma t\e^{2s}}{2+\alpha}\|\nabla u\|_2^2<4\pi$. Then, since $\|\cH(u,s)\|_\xi^\xi=\e^{(\xi-2)s}\| u\|_\xi^\xi$ for all $\xi\geq1$, from \eqref{estimate_convolution_term_H} one deduces
    $$\tJ(u,s)\les\frac{\e^{2s}}2\|\nabla u\|_2^2+\varepsilon\frac{\e^{2s}}2\|u\|^{4+\alpha}_{\frac{2(4+\alpha)}{2+\alpha}}+C_\varepsilon\e^{(2q-\frac{2+\alpha}{t'})s}\|u\|_{\frac{4qt'}{2+\alpha}}^{2q}\,,$$
    from which (i) follows, since the exponents are greater than $0$, provided $q$ is chosen large enough.

    In order to show (ii), for $v\in H^1(\R^2)$ let us define
    $$g(v):=\intd\left(I_\alpha\ast F(v)\right)F(v)\dd x$$
    and set $w(t):=g\left(\tfrac{tu}{\|u\|}\right)$, with $t=\e^s$. Then, ($f_3$) yields
    $$\frac{w'(t)}{w(t)}\geq\frac{2\mu}t\qquad\mbox{for}\quad t>0\,,$$
    from which
    $$g(tu)\geq g\left(\frac u{\|u\|}\right)\|u\|^{2\mu}t^{2\mu}.$$
    Therefore, one obtains
    $$J(\cH(u,s))\leq C_1\e^{2s}-C_2\e^{(2\mu-(2+\alpha))s}\to-\infty\qquad\mbox{as}\quad s\to+\infty\,,$$
    since $2\mu-(2+\alpha)>2$.
\end{proof}

\begin{Lem}\label{MP_geometry}
    Under ($f_0$)-($f_2$), there exists $K_a>0$ such that
    \begin{equation}\label{AB}
        0<\sup_{u\in A}J(u)<\inf_{u\in B}J(u)\,
    \end{equation}
    where
    \begin{equation*}
        A:=\left\{u\in\Sra\,\Big|\,\intd|\nabla u|^2\dd x\leq K_a\right\}\quad\mbox{and}\quad B:=\left\{u\in\Sra\,\Big|\,\intd|\nabla u|^2\dd x=2K_a\right\}.
    \end{equation*}
\end{Lem}
\begin{proof}
    Let $K_a<\tfrac{(2+\alpha)\pi}{2\gamma_0}$ and consider arbitrary $u_1\in A$ and $u_2\in B$. We aim at proving that $J(u_2)-J(u_1)>0$, from which \eqref{AB} easily follows. Starting from \eqref{estimate_convolution_term} and applying the Gagliardo-Nirenberg inequality \cite{Nirenberg1959}
    $$\|u\|_p\leq C_p\|\nabla u\|_2^{1-\frac2p}\|u\|_2^\frac2p,\qquad p\in(2,+\infty)\,,$$
    for $q,t>1$ and $\gamma>\gamma_0$ we infer
    \begin{equation*}
        \begin{split}
            \intd&\left(I_\alpha\ast F(u)\right)F(u)\dd x\les\varepsilon\|u\|_{\frac{2(\alpha+4)}{\alpha+2}}^{\alpha+4}+C_\varepsilon\left(\intd|u|^{\frac{4qt'}{2+\alpha}}\dd x\right)^\frac{2+\alpha}{2t'}\left(\intd\left(\e^{\frac{4\gamma t}{2+\alpha}u^2}-1\right)\dd x\right)^\frac{2+\alpha}{2t}\\
            &\les\varepsilon\|\nabla u\|_2^2\|u\|_2^{\alpha+2}+C_\varepsilon(q)\|\nabla u\|_2^{2q-\frac{2+\alpha}{t'}}\|u\|_2^\frac{2+\alpha}{t'}\left(\intd\left(\e^{\frac{4\gamma t}{2+\alpha}\|\nabla u\|_2^2\left(\frac{u}{\|\nabla u\|_2}\right)^2}-1\right)\dd x\right)^\frac{2+\alpha}{2t}.
        \end{split}
    \end{equation*}
    Since $u_2\in\Sra$ and $\|\nabla u_2\|_2^2=2K_a<\tfrac{(2+\alpha)\pi}{\gamma_0}$ small, then $\frac{4\gamma t}{2+\alpha}\|\nabla u_2\|_2^2<4\pi$ provided one chooses $t>1$ close to $1$ and $\gamma>\gamma_0$ close to $\gamma_0$, and therefore
    \begin{equation}\label{20240605e1}
        \intd\left(I_\alpha\ast F(u_2)\right)F(u_2)\dd x\leq\varepsilon C(a)\|\nabla u_2\|_2^2+C_\varepsilon(q,a)\|\nabla u_2\|_2^{2q-\frac{2+\alpha}{t'}}.
    \end{equation}
    Since, on the other hand, we have $\left(I_\alpha\ast F(u_1)\right)F(u_1)\geq0$, we deduce
    \begin{equation*}
        \begin{split}
            J(u_2)-J(u_1)&\geq\frac12\intd|\nabla u_2|^2\dd x-\frac12\intd|\nabla u_1|^2\dd x-\frac12\intd\left(I_\alpha\ast F(u_2)\right)F(u_2)\dd x\\
            &\geq\frac12 K_a-\varepsilon C_1K_a-C_2K_a^{q-\frac{2+\alpha}{2t'}}\\
            &\geq\frac14 K_a-C_2K_a^{q-\frac{2+\alpha}{2t'}}=:\rho_a>0\,,
        \end{split}
    \end{equation*}
    by taking $\varepsilon$ small, $q$ large, and $K_a$ small enough.  By arbitrariness of $u_1\in A$ and $u_2\in B$, we conclude that $\sup\limits_{u\in A}J(u)<\inf\limits_{u\in B}J(u)$. In particular, \eqref{20240605e1} is also valid after replacing $u_2$ by any $u_1\in A$, hence $J(u_1)>0$ for all $u_1\in A$.
\end{proof}

In view of Lemmas \ref{H_J} and \ref{MP_geometry}, we can now define the mountain pass level by
$$c_{mp}:=\inf_{h\in\Gamma}\sup_{t\in[0,1]}J(h(t))\,,$$
where $\Gamma=\left\{\gamma\in C\left([0,1],H^1_r(\R^2)\right)\,\left|\right.\,\gamma(0)\in A\,,\ J(\gamma(1))\leq0\right\}$.
\vskip0.2truecm
Since our problem is mass supercritical, it is well-known that mere Palais-Smale sequences are not necessary bounded in $H^1(\R^2)$. In order to retrieve a bounded sequence of almost critical points, following \cite{Jeanjean1997}, we look for a Palais-Smale-Poho\v zaev sequence, which is a Palais-Smale sequence with the additional property $\cP(u_n)\to0$, where $\cP$ is defined in \eqref{Pohozaev}. To this end, we need to prove that the auxiliary functional $\tJ:H^1(\R^2)\times\R\to\R$ defined in \eqref{tJ} has the same mountain pass structure as $J$, and that the two mountain pass levels coincide. Indeed, if so, similarly to the arguments in \cite[Proposition 2.2 and Lemma 2.4]{Jeanjean1997}, one can produce the existence of a (PSP)-sequence $(u_n,s_n)\in\Sra\times\R$ for $\tJ$, since $H^1_r(\R^2)$ is a Hilbert space, and derive from it a (PSP)-sequence for $J$ at the same level by taking $v_n:=\e^{s_n}u_n(\e^{s_n}\cdot)$. Following \cite[Lemma 4.1]{Dou2023}, and defining
$$\tc_{mp}=\inf_{\tgamma\in\tGamma}\sup_{t\in[0,1]}\tJ(\tgamma(t))\,,$$
where
$$\tGamma=\left\{\tgamma\in C\left([0,1],H^1_r(\R^2)\times\R\right)\,\left|\right.\,\tgamma(0)\in(A,0)\,,\ \tJ(\tgamma(1)))\leq0\right\}\,,$$
we prove the above claim in the following result:
\begin{Lem}
    Under ($f_0$)–($f_3$), the functional $\tJ$ has a mountain pass structure and $\tc_{mp}=c_{mp}\geq\rho_a>0\,$.
\end{Lem}
\begin{proof}
    If $\gamma\in\Gamma$, then $\tgamma:=(\gamma,1)\in\tGamma$ and $J(\gamma(t))=\tJ(\tgamma(1))$, so we directly infer $\tc_{mp}\leq c_{mp}$. On the other hand, for any $\tgamma(t)=(u(t),s(t))\in\tGamma$, then $u(0)\in A$, $s(0)=0$ and $\tJ(u(1),s(1)))\leq0$. Taking $\gamma(t):=\cH(u(t),s(t))\in\Sra$ for all $t\in[0,1]$ and $\gamma(0)=\cH(u(0),s(0))=u(0)\in A$ and $J(\gamma(1)))=J(\cH(u(1),s(1)))=\tJ(\tgamma(1))\leq 0\,.$
    As a result, $\gamma\in\Gamma$ and in turns $\tc_{mp}\geq c_{mp}$ holds.
\end{proof}
From the above reasons, the existence of a (PSP)-sequence for $J$ at level $c_{mp}$ follows.
\begin{Lem}\label{20240605L2}
    Under ($f_0$)–($f_3$), there exists a (PSP)-sequence $\{u_n\}\subset\Sra$, that is, such that
    $$J(u_n)\to c_{mp},\quad \left.J\right|'_{\Sra}(u_n)\to0,\quad\cP(u_n)\to0\quad\text{as}\quad n\to+\infty\,.$$
\end{Lem}
\noindent Here $\left.J\right|'_{\Sra}(u_n)\to0$ means $J'(u_n)[\varphi]=o_n(1)$ for all $\varphi\in H^1_r(\R^2)$ such that $\intd u_n\varphi\dd x=0\,$. Since $\{u_n\}\subset\Sra$, this implies that there exists a sequence $\{\lambda_n\}\subset\R$ such that
\begin{equation}\label{20240531-e2}
    \|J'(u_n)+\lambda_n N'(u_n)\|_{H^1_r(\R^2)'}=o_n(1)\,,
\end{equation}
where the functional $N:H_r^1(\R^2)\to\R$ is defined as
$$N(u):=\frac12\|u\|_2^2\,,$$
and $H^1_r(\R^2)'$ stands for the dual space of $H^1_r(\R^2)$.

\vskip0.2truecm
Since we are dealing with a critical exponential nonlinearity, in order to prove in the next section the uniform boundedness for the (PSP)-sequence found in Lemma \ref{20240605L2}, and avoid a bound from below like \eqref{global_bound_below_Schr}, we need a precise control on the mountain pass level in Lemma \ref{MP_level_bound} below. This may be accomplished by a delicate analysis on the Moser sequence given by
\begin{equation*}
    \tw_n(x)=\frac1{\sqrt{2\pi}}\begin{cases}
        \sqrt{\log n}\ \  &\mbox{for}\ \,0\leq|x|\leq\frac1n\,,\\
        \\
        \displaystyle \frac{\log(1/|x|)}{\sqrt{\log n}} & \mbox{for}\ \,\frac1n<|x|<1\,,\\
        \\
        0 & \mbox{for}\ \,|x|\geq1\,.
    \end{cases}
\end{equation*}
It is easy to see that
\begin{equation*}
    \|\tw_n\|_2^2=\frac1{4\log n}\left(1+o\left(\frac1{\log n}\right)\right)\qquad\mbox{and}\qquad\|\nabla\tw_n\|_2^2=1
\end{equation*}
for all $n\in\N$. Defining now $w_n:=a\frac{\tw_n}{\|\tw_n\|_2}$, so that $w_n\in\Sra$, for later use we compute
\begin{equation}\label{nabla_wn}
    \|\nabla w_n\|_2^2=\frac{a^2}{\|\tw_n\|_2^2}\intd|\nabla\tw_n|^2\dd x=4a^2\log n\left(1+o\left(\frac1{\log n}\right)\right)
\end{equation}
and, for $x\in B_\frac1n(0)$,
\begin{equation}\label{wn_ball}
    w_n(x)=\frac a{\sqrt{2\pi}}\frac{\sqrt{\log n}}{\|\tw_n\|_2}=a\sqrt{\frac2\pi}\log n\left(1+o\left(\frac1{\log n}\right)\right)\,.
\end{equation}
For a fixed $n\in\N$, let us introduce the function $g_n:\R^+\to\R$ as
\begin{align}\label{20240605e3}
    g_n(t)=\frac{t^2}2\|\nabla w_n\|_2^2-\frac1{2t^{2+\alpha}}\intd\left(I_\alpha\ast F(tw_n)\right)F(tw_n)\dd x\,,
\end{align}
and we claim that there exists $t_n>0$ such that $g_n(t_n)=\displaystyle{\max_{t>0}}\,g_n(t)$. Indeed, arguing similarly as in Lemma \ref{H_J} by means of \eqref{estF2}, one may show that
\begin{equation*}
    \begin{split}
	t^{-(2+\alpha)}\intd&\left(I_\alpha\ast F(tw_n)\right)F(tw_n)\dd x\les\varepsilon t^2\|w_n\|_{\frac{2(\alpha+4)}{(\alpha+2)}}^{\alpha+4}+C_\varepsilon t^{2q-\frac{2+\alpha}{t'}}\|w_n\|_{\frac{4qt'}{2+\alpha}}^{2q},
    \end{split}
\end{equation*}
which yields $g_n(t)\to-\infty$ as $t\to+\infty$, by choosing $q$ large enough and $g_n(t)\to0$ as $t\to0^+$. Since $g_n(t)>0$ for small $t$, the claim follows by continuity of $g_n$. We are now going to prove a suitable upper bound for $g_n(t_n)$.

\begin{Lem}\label{20240605-L1}
    Assuming ($f_0$)–($f_3$) and ($f_5$), there exists $n_0\in\N$ such that
    $$g_{n_0}(t_{n_0})=\displaystyle{\max_{t>0}}\,g_{n_0}(t)<\frac{(2+\alpha)\pi}{2\gamma_0}\,.$$
\end{Lem}
\begin{proof}
    Suppose by contradiction that
    \begin{equation}\label{MPlevel_contradiction}
	g_n(t_n)\geq\frac{(2+\alpha)\pi}{2\gamma_0}\,.
    \end{equation}
    for all $n\in\N$. Since $F\geq0$ by assumption, we have
    \begin{equation*}
	\frac{t_n^2}2\|\nabla w_n\|_2^2\geq g_n(t_n)\geq\frac{(2+\alpha)\pi}{2\gamma_0}\,,
    \end{equation*}
    and therefore, using \eqref{nabla_wn} we infer
    \begin{equation*}
	t_n^2\log n\geq\frac{(2+\alpha)\pi}{4a^2\gamma_0}\left(1+o\left(\frac1{\log n}\right)\right)\,.
    \end{equation*}
    Up to a subsequence, we may then assume that
    \begin{equation}\label{ell}
	\lim_{n\to+\infty}t_n^2\log n=\ell\in\left[\frac{(2+\alpha)\pi}{4a^2\gamma_0},+\infty\right)\cup\{+\infty\}\,.
    \end{equation}
    Since $t_n$ is a maximum point, from $t_ng_n'(t_n)=0$ and ($f_3$) we get
    \begin{equation}\label{tngntn_computations}
	\begin{split}
		t_n^2\|\nabla w_n\|_2^2&=\frac1{t_n^{2+\alpha}}\intd\!\left(I_\alpha\ast F(t_nw_n)\right)f(t_nw_n)t_nw_n\dd x-\frac{2+\alpha}{2t_n^{2+\alpha}}\intd\!\left(I_\alpha\ast F(t_nw_n)\right)F(t_nw_n)\dd x\\
		&\geq\left(\mu-\frac{2+\alpha}2\right)\frac1{t_n^{2+\alpha}}\intd\left(I_\alpha\ast F(t_nw_n)\right)F(t_nw_n)\dd x\\
		&\geq\left(\mu-\frac{2+\alpha}2\right)\frac1{t_n^{2+\alpha}}\int_{B_{\frac1n}(0)}\!\!\left(\int_{B_{\frac1n}(0)}\frac{c_\alpha F(t_nw_n(y))\dd y}{|x-y|^{2-\alpha}}\right)F(t_nw_n(x))\dd x\,.
	\end{split}
    \end{equation}
    Using \eqref{nabla_wn}, the fact that $w_n$ is constant in $B_\frac1n(0)$, see \eqref{wn_ball}, that for $x,y\in B_\frac1n(0)$ one has $\frac1{|x-y|^{2-\alpha}}\geq\left(\frac n2\right)^{2-\alpha}$, and taking into account the estimate from below at infinity for $F$ in \eqref{F_Ruf}, which holds by ($f_5$), we can derive
    \begin{equation*}
	\begin{split}
		4a^2t_n^2\log n&\left(1+o\left(\frac1{\log n}\right)\right)=t_n^2\|\nabla w_n\|_2^2\\
		&\geq\left(\mu-\frac{2+\alpha}2\right)\frac{c_\alpha}{t_n^{2+\alpha}}\frac{\pi^2}{n^4}\left(\frac n2\right)^{2-\alpha}\frac{(F(t_nw_n)(t_nw_n)^2)^2}{(t_nw_n)^4}\Bigg|_{w_n=a\sqrt{\tfrac2\pi}\log n\left(1+o\left(\frac1{\log n}\right)\right)}\\
		  &\geq\frac{\beta_0^2c_\alpha\pi^2}{2^{4-\alpha}\gamma_0^2}\left(\mu-\frac{2+\alpha}2\right)\frac{n^{-(2+\alpha)}}{t_n^{6+\alpha}}\frac{\exp\left\{4\gamma_0\frac{a^2}\pi t_n^2\log^2n\left(1+o\left(\frac1{\log n}\right)\right)\right\}}{\frac{4a^4}{\pi^2}\log^4n\left(1+o\left(\frac1{\log n}\right)\right)}\,,
	\end{split}
    \end{equation*}
    from which
    \begin{equation*}
	\begin{split}
		4a^2&\geq\frac{K_1\exp\left\{\left[K_2\left(t_n^2\log n\right)\left(1+o\left(\frac1{\log n}\right)\right)-(2+\alpha)\right]\log n\right\}}{\left(t_n^2\log n\right)^{4+\frac\alpha2}\left(\log n\right)^{1-\frac\alpha2}\left(1+o\left(\frac1{\log n}\right)\right)}\,,
	\end{split}
    \end{equation*}
    where
    $$K_1:=\frac{\beta_0^2c_\alpha\pi^4}{2^{6-\alpha}\gamma_0^2a^4}\left(\mu-\frac{2+\alpha}2\right)\quad\ \mbox{and}\quad\  K_2:=4\gamma_0\frac{a^2}\pi\,.$$
    It is then easy to rule out the possibility in \eqref{ell} that $\ell=+\infty$. Hence $\ell<+\infty$ and
    \begin{equation*}
	4a^2\geq\frac{K_1\exp\left\{\left[K_2(\ell+o_n(1))\left(1+o\left(\frac1{\log n}\right)\right)-(2+\alpha)\right]\log n\right\}}{(\ell+o_n(1))^{4+\frac\alpha2}\left(\log n\right)^{1-\frac\alpha2}\left(1+o\left(\frac1{\log n}\right)\right)}\,,
    \end{equation*}
    which is again a contradiction unless $K_2\ell\leq2+\alpha$, that is, by \eqref{ell},
    \begin{equation*}
	\lim_{n\to+\infty}t_n^2\log n=\ell=\frac{(2+\alpha)\pi}{4a^2\gamma_0}\,.
    \end{equation*}
    Similarly as in \eqref{tngntn_computations}, we can also obtain
    \begin{equation*}
	\begin{split}
		\frac1{t_n^{2+\alpha}}\intd&\left(I_\alpha\ast F(t_nw_n)\right)F(t_nw_n)\dd x\geq\frac{c_\alpha\pi^2}{t_n^{2+\alpha}n^4}\left(\frac n2\right)^{2-\alpha}\left(F(t_nw_n)\right)^2\\
		&\geq\frac{c_\alpha\pi^2\beta_0^2}{2^{4-\alpha}\gamma_0^2t_n^{2+\alpha}}\exp\left\{4\gamma_0\frac{a^2}{\pi}t_n^2\log^2n\left(1+o\left(\frac1{\log n}\right)\right)\right\}\frac{n^{-(2+\alpha)}}{(t_nw_n)^4}\\
		&=\frac{c_\alpha\pi^4\beta_0^2}{a^42^{6-\alpha}\gamma_0^2}\frac1{t_n^{6+\alpha}\log^4n}\exp\left\{\left[4\gamma_0\frac{a^2}{\pi}t_n^2\log n\left(1+o\left(\frac1{\log n}\right)\right)-(2+\alpha)\right]\log n\right\}\\
		  &\geq\frac{2\tK_1}{\left(\log n\right)^{1-\frac\alpha2}}\exp\left\{\left[\tK_2^{(n)}t_n^2\log n-(2+\alpha)\right]\log n\right\}\,,
	\end{split}
    \end{equation*}
    where
    $$\tK_1:=\frac{c_\alpha\pi^4\beta_0^2\left(\ell^{-\left(3+\frac\alpha2\right)}-\varepsilon\right)}{a^42^{7-\alpha}\gamma_0^2}\qquad\mbox{and}\qquad\tK_2^{(n)}:=\frac{4\gamma_0a^2}{\pi}\left(1+o\left(\frac1{\log n}\right)\right)$$
    with $\varepsilon>0$ small and $n$ large enough. Hence, by \eqref{nabla_wn} and the definition of $g_n$ in \eqref{20240605e3}, we obtain
    \begin{equation*}
	\begin{split}
		g_n(t_n)&\leq2a^2t_n^2\log n\left(1+o\left(\frac1{\log n}\right)\right)-\frac{\tK_1}{\left(\log n\right)^{1-\frac\alpha2}}\exp\left\{\left[\tK_2^{(n)}t_n^2\log n-(2+\alpha)\right]\log n\right\}\\
		&=:\tK_0^{(n)}s-\frac{\tK_1}{\left(\log n\right)^{1-\frac\alpha2}}\exp\{\left[\tK_2^{(n)}s-(2+\alpha)\right]\log n\}=:h_n(s)\,
	\end{split}
    \end{equation*}
    having set $s:=t_n^2\log n$ and $\tK_0^{(n)}:=2a^2\left(1+o\left(\tfrac1{\log n}\right)\right)$, and we are going to maximise the function $h_n$. Defining $s_n:=\displaystyle{\max_{s>0}h_n(s)}$, we easily infer that $s_n$ fulfils
    \begin{equation}\label{s_n}
	n^{\tK_2^{(n)}s_n-(2+\alpha)}\left(\log n\right)^\frac\alpha2=\frac{\tK_0^{(n)}}{\tK_1\tK_2^{(n)}}\,.
    \end{equation}
    Hence necessarily $\{s_n\}_n$ is bounded and we can assume, up to a subsequence, that $s_n\to s^*\in[0,+\infty)$. Our goal is now to show that
    \begin{equation*}
	s^*=\ell=\frac{(2+\alpha)\pi}{4\gamma_0a^2}\,.
    \end{equation*}
    Suppose by contradiction that $s^*<\frac{(2+\alpha)\pi}{4\gamma_0a^2}$, then
    $$\tK_2^{(n)}s_n-(2+\alpha)=\frac{4a^2\gamma_0}\pi\left(1+o\left(\frac1{\log n}\right)\right)s_n-(2+\alpha)<-\varepsilon_0<0$$
    for large $n$, and this would be in contradiction with \eqref{s_n}. The case $s^*>\frac{(2+\alpha)\pi}{4\gamma_0a^2}$ can be treated analogously and ends up in a similar contradiction, thus the claim is proved, and we can write $s_n=\ell+s_n^*$ with $s_n^*=o_n(1)$. Therefore,
    $$\tK_2^{(n)}s^*-(2+\alpha)=(2+\alpha)o\left(\frac1{\log n}\right)\,,$$
    and then, by \eqref{s_n},
    $$\frac{\tK_0^{(n)}}{\tK_1\tK_2^{(n)}}=n^{\left[\tK_2^{(n)}s^*-(2+\alpha)\right]+K_2s_n^*}\left(\log n\right)^\frac\alpha2=n^{\tK_2^{(n)}s_n^*+o\left(\frac1{\log n}\right)}\left(\log n\right)^\frac\alpha2.$$
    This leads to an other contradiction unless $s_n^*=o\left(\frac1{\log n}\right)$. But then, again from \eqref{s_n}
    \begin{equation*}
	\begin{split}
		g_n(t_n)&\leq h_n(t_n^2\log n)\leq h_n(s_n)=\tK_0^{(n)}s_n-\frac{\tK_0^{(n)}}{\tK_2^{(n)}\log n}\\
		&=\tK_0^{(n)}\left[\frac{(2+\alpha)\pi}{4\gamma_0a^2}+o\left(\frac1{\log n}\right)-\frac1{\tK_2^{(n)}\log n}\right]\\
		&<\tK_0^{(n)}\frac{(2+\alpha)\pi}{4\gamma_0a^2}=\frac{(2+\alpha)\pi}{2\gamma_0}\,,
	\end{split}
    \end{equation*}
    which finally contradicts \eqref{MPlevel_contradiction} and concludes the proof.
\end{proof}

\begin{Cor}\label{MP_level_bound}
    Under ($f_0$)–($f_3$) and ($f_5$), we have $c_{mp}<\frac{(2+\alpha)\pi}{2\gamma_0}\,.$
\end{Cor}
\begin{proof}
    For $n_0\in\N$ set as in Lemma \ref{20240605-L1}, by Lemma \ref{H_J} there exist $s_1<-1$ sufficiently negative and $s_2>1$ large enough such that $\cH(w_{n_0},s_1)\in A$ and $J(\cH(w_{n_0},s_2))<0.$ Let $h_0(t):=\cH(w_{n_0},ts_2+(1-t)s_1)$ with $t\in[0,1]$. Then $h_0\in\Gamma$, and thus
    \begin{align*}
        c_{mp}\leq\max_{t\in [0,1]}J(h_0(t))=\max_{s_1\leq s\leq s_2}J(\cH(w_{n_0},s))\leq \max_{t\geq0}g_{n_0}(t)<\frac{(2+\alpha)\pi}{2\gamma_0}\,.
    \end{align*}
\end{proof}

\section{The $(PSP)_{c_{mp}}-$ condition}\label{sec_PSP}
In this section, we verify that the $(PSP)_{c_{mp}}$ - condition for the constrained functional $J\big|_{\Sra}$ holds. We first establish some convergence results.

\begin{Lem}\label{lemma 2.3}
    Assume that $\{u_n\}\subset H_r^1(\R^2)$ such that $u_n\rightharpoonup u_0$ weakly in $H_r^1(\R^2)$ and
    \begin{equation}\label{Assumption_Lemma4.1}
        \limsup_{n\to+\infty}\|\nabla(u_n-u_0)\|_2^2<\frac{(2+\alpha)\pi}{\gamma_0}\,.
    \end{equation}
    Under the assumptions ($f_0$)-($f_3$) and ($f_5$), there holds
    \begin{align}\label{20240531-A4}
        \int_{\R^2}(I_\alpha\ast F(u_n))F(u_n)\dd x=\int_{\R^2}(I_\alpha\ast F(u_0))F(u_0)\dd x+o_n(1)\,,
    \end{align}
    and
    \begin{align}\label{20240531-A5}
        \int_{\R^2}(I_\alpha\ast F(u_n))f(u_n)u_n\dd x=\int_{\R^2}(I_\alpha\ast F(u_0))f(u_0)u_0\dd x+o_{n}(1)\,.
    \end{align}
\end{Lem}
\begin{proof}
    Recalling \eqref{estf1}, by Proposition \ref{Cao_ineq} (i), one can show that $f(u_0)\in L^\tau(\R^2)$ for $\tau>1$ close to $1$. Moreover, since $H^1_r(\R^2)$ is compactly embedded in $L^q(\R^2)$ for $q>2$ by \cite{Lions82}, we have that $u_n\to u_0$ in $L^{\tau^\prime}$ with $\tau^{\prime}=\frac{\tau}{\tau-1}>2$. Thus,
    \begin{align}\label{20240531-A1}
        \int_{\R^2}|f(u_0)(u_n-u_0)|\dd x\leq\|f(u_0)\|_\tau\|u_n-u_0\|_{\tau^{\prime}}=o_n(1)\,.
    \end{align}
	By \cite[Lemma 4.1]{ACTY}, we know
	\begin{align}\label{20240531-A3}
		\|I_\alpha\ast F(u_0)\|_{\infty}<+\infty\,.
	\end{align}
	We note that, although our assumptions on $f$ are weaker than those in \cite{ACTY}, an upper bound on the nonlinearity as in \eqref{estf1} is sufficient to proceed as in the proof of \cite[Lemma 4.1]{ACTY} and obtain \eqref{20240531-A3}. By combining \eqref{20240531-A3} with \eqref{estf1}, \eqref{Assumption_Lemma4.1}, and Lemma \ref{lemma 2.2}, one infers that, up to a subsequence, $(I_\alpha\ast F(u_0))f(u_n)\rightharpoonup (I_\alpha\ast F(u_0))f(u_0)$ weakly in $L^{\tau}(\R^2)$ for $\tau>1$ close to $1$. Therefore,
    \begin{align}\label{20240531-A2}
        \int_{\R^2}(I_\alpha\ast F(u_0))(f(u_n)-f(u_0))u_n\dd x=o_n(1)\,.
    \end{align}
     As a result, by \eqref{20240531-A3} and the Hardy-Littlewood-Sobolev inequality, one has
    \begin{equation}\label{eq_number}
    	\begin{split}
    		&\left|\int_{\R^2}(I_\alpha\ast F(u_n))f(u_n)u_n\dd x-\int_{\R^2}(I_\alpha\ast F(u_0))f(u_0)u_0\dd x\right|\\
    		&\quad\leq\left|\int_{\R^2}(I_\alpha\ast (F(u_n)-F(u_0)))f(u_n)u_n\dd x\right|+\left|\int_{\R^2}(I_\alpha\ast F(u_0))(f(u_n)-f(u_0))u_n\dd x\right|\\
    		&\qquad+\left|\int_{\R^2}(I_\alpha\ast F(u_0))f(u_0)(u_n-u_0)\dd x\right|\\
    		&\quad\les\|F(u_n)-F(u_0)\|_{\frac4{2+\alpha}}\|f(u_n)u_n\|_{\frac4{2+\alpha}}+\left|\int_{\R^2}(I_\alpha\ast F(u_0))(f(u_n)-f(u_0))u_n\dd x\right|\\
    		&\qquad+\int_{\R^2}|f(u_0)||u_n-u_0|\dd x\\
    		&\quad=\|F(u_n)-F(u_0)\|_{\frac4{2+\alpha}}\|f(u_n)u_n\|_{\frac4{2+\alpha}}+o_n(1)\,,
    	\end{split}
    \end{equation}
	thanks to \eqref{20240531-A2} and \eqref{20240531-A1}. By combining \eqref{estf1}, the continuous embedding $H^1_r(\R^2)\hookrightarrow L^q(\R^2)$ for all $q>2$, and the uniform bound given by Lemma \ref{lemma 2.2} with $u_0=0$, one obtains
    \begin{equation}\label{funun}
        \|f(u_n)u_n\|_{\frac4{2+\alpha}}\leq C\,.
    \end{equation}
    On the other hand, by the mean value theorem, the H\"older inequality, and ($f_1$)–($f_2$), there exists a function $\theta:\R^2\to[0,1]$ such that
    \begin{equation*}
        \begin{split}
            \|F(u_n)-F(u_0)&\|_{\frac4{2+\alpha}}^{\frac4{2+\alpha}}=\int_{\R^2}\big|f(u_0+\theta(u_n-u_0))(u_n-u_0)\big|^{\frac4{2+\alpha}}\dd x\\
            &\les\int_{\R^2}|u_0+\theta(u_n-u_0)|^2|u_n-u_0|^{\frac4{2+\alpha}}\dd x\\
            &~~~~+C_\varepsilon\int_{\R^2}|u_0+\theta(u_n-u_0)|^{\frac{4(q-1)}{2+\alpha}}|u_n-u_0|^{\frac{4}{2+\alpha}}(\e^{\gamma|u_0+\theta(u_n-u_0)|^2}-1)^{\frac{4}{2+\alpha}}\dd x\\
            &\leq o_n(1)+\Big[\int_{\R^2}|u_0+\theta(u_n-u_0)|^{\frac{4\tau'
            (q-1)}{2+\alpha}}|u_n-u_0|^{\frac{4\tau'}{2+\alpha}}\dd x\Big]^{\frac{1}{\tau'}}\\
            &\quad\times\Big[\int_{\R^2}(\e^{\gamma|u_0+\theta(u_n-u_0)|^2}-1)^{\frac{4\tau}{2+\alpha}}\dd x\Big]^{\frac1\tau}\\
            &\leq o_n(1)+\left(\intd\left(|u_0|+|u_n|\right)^{\frac{4\tau'\nu'(q-1)}{2+\alpha}}\dd x\right)^{\frac1{\tau'\nu'}}\left(\intd|u_n-u_0|^{\frac{4\tau'\nu}{2+\alpha}}\dd x\right)^{\frac1{\tau'\nu}}\\
            &\quad\times\Big[\int_{\R^2}(\e^{\gamma|u_0+\theta(u_n-u_0)|^2}-1)^{\frac{4\tau}{2+\alpha}}\dd x\Big]^{\frac1\tau},
        \end{split}
    \end{equation*}
    where $\frac1\nu+\frac1{\nu'}=1$. We note that the first factor is bounded by weak convergence, the second is $o_n(1)$ because of the compact embedding. Moreover, by using the identity \eqref{exp_identity} and the H\"older inequality, one can also bound the last term as
    \begin{equation*}
        \begin{split}
            \intd(\e^{\gamma|u_0+\theta(u_n-u_0)|^2}&-1)^{\frac{4\tau}{2+\alpha}}\dd x\leq\intd\Big(\e^{\frac{4\tau\gamma}{2+\alpha}|u_0+\theta(u_n-u_0)|^2}-1\Big)\dd x\\
            &\leq\intd\Big(\e^{\frac{4\gamma\tau(1+\frac1\sigma)}{2+\alpha}u_0^2}-1\Big)\Big(\e^{\frac{4\gamma\tau(1+\sigma)}{2+\alpha}|u_n-u_0|^2}-1\Big)\dd x\\
            &\quad+\intd\Big(\e^{\frac{4\gamma\tau(1+\frac1\sigma)}{2+\alpha}u_0^2}-1\Big)\dd x+\intd\Big(\e^{\frac{4\gamma\tau(1+\sigma)}{2+\alpha}|u_n-u_0|^2}-1\Big)\dd x\\
            &\leq\left(\intd\Big(\e^{\frac{4\gamma\tau(1+\frac1\sigma)\eta'}{2+\alpha}u_0^2}-1\Big)\dd x\right)^\frac1{\eta'}\left(\intd\Big(\e^{\frac{4\gamma\tau(1+\sigma)\eta}{2+\alpha}|u_n-u_0|^2}-1\Big)\dd x\right)^\frac1\eta\\
            &\quad+\intd\Big(\e^{\frac{4\gamma\tau(1+\frac1\sigma)}{2+\alpha}u_0^2}-1\Big)\dd x+\intd\Big(\e^{\frac{4\gamma\tau(1+\sigma)}{2+\alpha}|u_n-u_0|^2}-1\Big)\dd x\,.
        \end{split}
    \end{equation*}
    Noticing that Lemma \ref{lemma 2.2} yields
    $$\intd\left(\e^{\frac{4\gamma\tau(1+\sigma)}{2+\alpha}|u_n-u_0|^2}-1\right)\dd x\leq C\qquad\mbox{and}\qquad\intd\Big(\e^{\frac{4\gamma\tau(1+\sigma)\eta}{2+\alpha}|u_n-u_0|^2}-1\Big)\dd x\leq C\,,$$
    provided $\tau>1$ and $\eta>1$ are both close to $1$, and $\sigma>0$ is close to $0$, while Proposition \ref{Cao_ineq} (i) implies
    $$\intd\left(\e^{\frac{4\gamma\tau(1+\frac1\sigma)}{2+\alpha}u_0^2}-1\right)\dd x<+\infty\qquad\mbox{and}\qquad\intd\Big(\e^{\frac{4\gamma\tau(1+\frac1\sigma)\eta'}{2+\alpha}u_0^2}-1\Big)\dd x<+\infty\,,$$
    one then infers
    \begin{equation}\label{20240531-A6}
        \|F(u_n)-F(u_0)\|_{\frac4{2+\alpha}}=o_n(1)\,.
    \end{equation}
    As a result, \eqref{20240531-A5} follows by combining \eqref{eq_number} with \eqref{funun} and \eqref{20240531-A6}. Using an argument similar to \cite[Corollary 2.5]{Dou2023}, one has
    \begin{align}\label{20240531-A7}
        \int_{\R^2}|F(u_n)-F(u_0)|\dd x=o_n(1)\,.
    \end{align}
    Hence, using the Hardy-Littlewood-Sobolev inequality and \eqref{20240531-A3} again, we derive
    \begin{align*}
        &\left|\int_{\R^2}(I_\alpha\ast F(u_n))F(u_n)\dd x-\int_{\R^2}(I_\alpha\ast F(u_0))F(u_0)\dd x\right|\\
        &\leq\left|\int_{\R^2}(I_\alpha\ast (F(u_n)-F(u_0)))F(u_n)\dd x\right|+\left|\int_{\R^2}(I_\alpha\ast F(u_0))(F(u_n)-F(u_0))\dd x\right|\\
        &\les\|F(u_n)-F(u_0)\|_{\frac4{2+\alpha}}\|F(u_n)\|_{\frac4{2+\alpha}}+\int_{\R^2}|F(u_n)-F(u_0)|\dd x\,.
    \end{align*}
    Since $\|F(u_n)\|_{\frac4{2+\alpha}}\leq C$ by arguing similarly to \eqref{funun}, by \eqref{20240531-A6} and \eqref{20240531-A7}, then \eqref{20240531-A4} holds.
\end{proof}

We are now ready to establish the $(PSP)_{c_{mp}}-$ condition.
\begin{Lem}\label{lemma 3.6}
    Suppose $\{u_n\}\subset\Sra$ is a $(PSP)_{c_{mp}}-$ sequence with $c_{mp}<\frac{(2+\alpha)\pi}{2\gamma_0}$, and $c_{mp}\neq0,$ i.e.,
    \begin{equation}\label{un_PSP}
        J(u_n)\to c_{mp}<\frac{(2+\alpha)\pi}{2\gamma_0},~~~J\big|_{\Sra}^{\prime}(u_n)\to 0~~\text{and}~~\cP(u_n)\to 0 ~\hbox{as}~n\to+\infty\,.
    \end{equation}
    Assuming ($f_0$)–($f_5$), then $\{u_n\}$ is bounded with
    \begin{equation}\label{bound_un}
        \|\nabla u_n\|_2^2\leq\frac{2c_{mp}(\mu-\frac{2+\alpha}2)}{\mu-(2+\frac\alpha2)}+o_n(1)
    \end{equation}
    and, up to a subsequence, $\{u_n\}$ converges strongly in $H_r^1(\R^2)$ to $u_a$, which is a weak solution of \eqref{eq} for some $\lambda_a>0$.
\end{Lem}
\begin{proof}
	By \eqref{un_PSP} and ($f_3$) one has
    \begin{equation*}
        \begin{split}
            o_n(1)&=\cP(u_n)=\|\nabla u_n\|_2^2+\frac{2+\alpha}2\intd(I_\alpha\ast F(u_n))F(u_n)\dd x-\intd(I_\alpha\ast F(u_n))f(u_n)u_n\dd x\\
            &\leq\|\nabla u_n\|_2^2-\left(\mu-\frac{2+\alpha}2\right)\intd(I_\alpha\ast F(u_n))F(u_n)\dd x
        \end{split}
    \end{equation*}
    and so
    \begin{equation*}
        \intd(I_\alpha\ast F(u_n))F(u_n)\dd x\leq\frac{\|\nabla u_n\|_2^2}{\mu-\frac{2+\alpha}2}+o_n(1)\,.
    \end{equation*}
    Hence, using $J(u_n)\to c_{mp}$, one infers
    $$
        c_{mp}\geq\frac12\|\nabla u_n\|_2^2-\frac1{2\mu-(2+\alpha)}\|\nabla u_n\|_2^2+o_n(1)=\frac{\mu-(2+\frac\alpha2)}{2(\mu-\frac{2+\alpha}{2})}\|\nabla u_n\|_2^2+o_n(1)\,.
    $$
    Since $\mu>2+\frac{\alpha}2$, the sequence $\{\|\nabla u_n\|_2^2\}$ is bounded and \eqref{bound_un} holds. This easily implies that also $\{\|u_n\|_{H^1(\R^2)}\}$ is bounded since $\{u_n\}\subset\Sra$. Moreover, we claim that
    \begin{equation}\label{bdd}
        \int_{\R^2}(I_\alpha\ast F(u_n))F(u_n)\dd x\leq C\,,\quad\int_{\R^2}(I_\alpha\ast F(u_n))f(u_n)u_n\dd x\leq C\,,\quad\mbox{and}\quad0\leq\lambda_n\leq C\,.
    \end{equation}
    Indeed, the first is a consequence of \eqref{bound_un} and $J(u_n)\to c_{mp}$; the second is implied by the first and $\cP(u_n)\to0$; the third follows, since
    \begin{equation}\label{lambda_n}
        \begin{split}
            \lambda_n&=-\frac{J'(u_n)u_n}{N'(u_n)u_n}=-\frac{\|\nabla u_n\|_2^2-\int_{\R^2}(I_\alpha\ast F(u_n))f(u_n)u_n\dd x}{a^2}\\
            &=\frac1{a^2}\left(1+\frac\alpha2\right)\int_{\R^2}(I_\alpha\ast F(u_n))F(u_n)\dd x+o_n(1)\,,
        \end{split}
    \end{equation}
    again by $\cP(u_n)\to0$. Hence, going to a subsequence, we may assume that $u_n\rightharpoonup u_a$ weakly in $H^1(\R^2)$ for some $u_a\in H_r^1(\R^2)$, and $\lambda_n\to\lambda_a$ as $n\to+\infty$.

    We firstly prove that $u_a\neq 0$. If not, namely supposing $u_a\equiv0$, by the radial compact embedding we infer $\|u_n\|_{2+\frac{\alpha}2}=o_n(1)$. By \eqref{f1f4}, for any $\varepsilon>0$, we can find $C_\varepsilon>0$ such that
    \begin{align}\label{20240531-e3}
        \int_{\R^2}(I_\alpha\ast F(u_n))F(u_n)\dd x &\leq \varepsilon\int_{\R^2}(I_\alpha\ast F(u_n))f(u_n)u_n\dd x+C_\varepsilon\|u_n\|_{2+\frac{\alpha}{2}}^{2+\frac\alpha2}\nonumber\\
        &\leq\varepsilon\int_{\R^2}(I_\alpha\ast F(u_n))f(u_n)u_n\dd x+o_n(1)\,.
    \end{align}
    Since $\{\int_{\R^2}(I_\alpha\ast F(u_n))f(u_n)u_n\dd x\}$ is bounded, by the arbitrariness of $\varepsilon$, we conclude that
    \begin{align*}
        \lim_{n\to+\infty}\int_{\R^2}(I_\alpha\ast F(u_n))F(u_n)\dd x=0\,.
    \end{align*}
    By \eqref{un_PSP}, this yields
    \begin{align}\label{20240531-e3-1}
        c_{mp}=\frac12\lim_{n\to+\infty}\|\nabla u_n\|_2^2
    \end{align}
    and
    \begin{align*}
        \lim_{n\to+\infty}\|\nabla u_n\|_2^2<\frac{(2+\alpha)\pi}{\gamma_0}\,.
    \end{align*}
    Note that, by Lemma \ref{lemma 2.3}, we also have
    \begin{align*}
        \lim_{n\to+\infty}\intd(I_\alpha\ast F(u_n))f(u_n)u_n\dd x=0\,.
    \end{align*}
    Combining the above with $\cP(u_n)\to 0$, we obtain that $\|\nabla u_n\|_2^2\to 0$, which, by \eqref{20240531-e3-1} is a contradiction to the condition $c_{mp}\neq 0$. Hence $u_a\not\equiv 0$ is proved. Since $u_n\to u_a$ in $L^{2+\frac\alpha2}(\R^2)$ and that $\{(I_\alpha\ast F(u_n))f(u_n)u_n\}$ is bounded in $L^1(\R^2)$, then $\int_{\R^2}(I_\alpha\ast F(u_n))F(u_n)\dd x$ is uniformly bounded and moreover possesses the uniform integrability condition. Indeed, fixed $\widetilde\varepsilon>0$, if $\delta>0$ and $X\subset\R^2$ so that $|X|<\delta$, then by \eqref{20240531-e3} and H\"older's inequality with $\tau,\tau'>1$, one has
    \begin{equation}\label{unif_int_cond_X}
        \begin{split}
            \int_X(I_\alpha\ast F(u_n))F(u_n)\dd x &\leq \varepsilon\int_X(I_\alpha\ast F(u_n))f(u_n)u_n\dd x+C_\varepsilon\int_X|u_n|^{2+\frac\alpha2}\dd x\\
            &\leq\varepsilon\intd(I_\alpha\ast F(u_n))f(u_n)u_n\dd x+|X|^\frac1{\tau'}\left(\intd|u_n|^{\left(2+\frac\alpha2\right)\tau}\dd x\right)^\frac1\tau\\
            &\leq C\varepsilon+C\delta^\frac1{\tau'}<\widetilde\varepsilon\,,
        \end{split}
    \end{equation}
    provided $n$ large enough and $\varepsilon$ and $\delta$ small enough. Similarly, for $R>0$ one has
    \begin{equation}\label{tightness}
        \begin{split}
            \int_{B_R(0)^c}(I_\alpha\ast F(u_n))F(u_n)\dd x &\leq \varepsilon\intd(I_\alpha\ast F(u_n))f(u_n)u_n\dd x+C_\varepsilon\int_{B_R(0)^c}|u_n|^{2+\frac\alpha2}\dd x\\
            &\leq C\varepsilon+C_\varepsilon\left(\int_{B_R(0)^c}|u_0|^{2+\frac\alpha2}\dd x+\intd|u_n-u_0|^{2+\frac\alpha2}\dd x\right).
        \end{split}
    \end{equation}
    Fixing $\widetilde\varepsilon>0$, and in turn $\varepsilon<\frac{\widetilde\varepsilon}{2C}$ one can find $n$ large and $R$ large so that 
    $$\int_{B_R(0)^c}|u_0|^{2+\frac\alpha2}\dd x+\intd|u_n-u_0|^{2+\frac\alpha2}\dd x<\frac{\widetilde\varepsilon}{2C_\varepsilon}\,,$$
    from which one deduces also the tightness property of $\{\int_{\R^2}(I_\alpha\ast F(u_n))F(u_n)\dd x\}$. As a consequence of \eqref{unif_int_cond_X} and \eqref{tightness} one may apply Vitali's convergence theorem (see \cite[Chapter 5]{RF}), and obtain
    \begin{equation}\label{20240531-e4}
        \begin{split}
            \lim_{n\to+\infty}\int_{\R^2}(I_\alpha\ast F(u_n))F(u_n)\dd x&=\int_{\R^2}\lim_{n\to+\infty}(I_\alpha\ast F(u_n))F(u_n)\dd x\\
            &=\int_{\R^2}(I_\alpha\ast F(u_a))F(u_a)\dd x\,.
        \end{split}
    \end{equation}
	By \eqref{lambda_n}, this yields
    \begin{equation}\label{20240531-e7}
        \begin{split}
            \lambda_a&=\lim_{n\to+\infty}\lambda_n=\frac1{a^2}\left(1+\frac\alpha2\right)\int_{\R^2}(I_\alpha\ast F(u_a))F(u_a)\dd x>0\,.
        \end{split}
    \end{equation}
    We show now that $u_a$ is a weak solution of
    \begin{align}\label{20240531-e6}
        -\Delta u_a+\lambda_au_a=(I_\alpha\ast F(u_a))f(u_a)\quad\ \text{ in }\ \R^2
    \end{align}
    by following the line of \cite[Lemma 2.4]{ACTY}. In fact, for $\varphi\in C^\infty_0(\R^2)$, defining $v_n:=\frac\varphi{1+|u_n|}$, we have
    \begin{equation*}
        \begin{split}
            \|v_n\|_{H^1(\R^2)}^2&=\intd\frac{\varphi^2}{(1+|u_n|)^2}\dd x+\intd\left|\frac{\nabla\varphi}{1+|u_n|}-\frac{\varphi|\nabla u_n|}{(1+|u_n|)^2}\right|^2\dd x\\
            &\leq\intd\left(\frac{\varphi^2+|\nabla\varphi|^2}{(1+|u_n|)^2}+2\frac{|\nabla\varphi\nabla u_n|}{(1+|u_n|)^3}+\frac{|\varphi|^2|\nabla u_n|^2}{(1+|u_n|)^4}\right)\dd x\\
            &\leq\|\varphi\|_2^2+2\|\nabla \varphi\|_2^2+(1+\|\varphi\|_\infty^2)\|\nabla u_n\|_2^2\,,
        \end{split}
    \end{equation*}
    hence $v_n\in H^1_r(\R^2)$ for all $n\in\N$ and can be chosen as a test function in \eqref{20240531-e2}. Moreover, by \eqref{bound_un}, we obtain that $\{v_n\}$ is uniformly bounded in $H^1_r(\R^2)$. Let now $\Omega\subset\subset\R^2$ and $\varphi\in C^\infty_0(\R^2)$ so it takes values in $[0,1]$ and $\varphi\equiv1$ on $\Omega$. We have
    \begin{equation}\label{Zoldo}
        \begin{split}
            \int_\Omega\left(I_\alpha\ast F(u_n)\right)f(u_n)\dd x&\leq\int_{\Omega\cap\{|u_n|>1\}}\left(I_\alpha\ast F(u_n)\right)f(u_n)|u_n|\dd x+2\int_{\Omega\cap\{|u_n|\leq 1\}}\left(I_\alpha\ast F(u_n)\right)\frac{f(u_n)}{1+|u_n|}\dd x\\
            &\leq\intd\left(I_\alpha\ast F(u_n)\right)f(u_n)u_n\dd x+2\intd\left(I_\alpha\ast F(u_n)\right)f(u_n)\frac{\varphi}{1+|u_n|}\dd x\,,
        \end{split}
    \end{equation}
    where ($f_0$) has been used in the first term. Combining \eqref{bdd} with
    \begin{equation*}
        \begin{split}
            \intd\left(I_\alpha\ast F(u_n)\right)f(u_n)\frac{\varphi}{1+|u_n|}\dd x&=\intd\nabla u_n\nabla\left(\frac\varphi{1+|u_n|}\right)\dd x+\lambda_n\intd \frac{u_n\varphi}{1+|u_n|}\dd x+o\left(\|v_n\|_{H^1(\R^2)}\right)\\
            &\leq\intd\left(\frac{|\nabla u_n\nabla\varphi|}{1+|u_n|}+\frac{|\nabla u_n|^2\varphi}{(1+|u_n|)^2}\right)\dd x+(\lambda_a+1)\intd\frac{|u_n|\varphi}{1+|u_n|}\dd x+o_n(1)\\
            &\les\|u_n\|_{H^1(\R^2)}^2+\|\varphi\|_{H^1(\R^2)}^2\leq C\,,
        \end{split}
    \end{equation*}
    from \eqref{Zoldo} we infer that the measure $\nu_n$ defined by
    $$\nu_n(\Omega):=\int_\Omega\left(I_\alpha\ast F(u_n)\right)f(u_n)\dd x$$
    has uniformly bounded total variation, hence there exists a measure $\nu$ such that, up to a subsequence, $\nu_n\stackrel{*}{\rightharpoonup}\nu$, namely
    $$\int_\Omega\left(I_\alpha\ast F(u_n)\right)f(u_n)\,\varphi\dd x\to\int_\Omega\varphi\dd\nu$$
    for all $\varphi\in C^\infty_0(\Omega)$. As in \cite[Lemma 2.4]{ACTY} we may then conclude that $\nu$ is absolutely continuous with respect to the Lebesgue measure and it can be identified as $\nu=\left(I_\alpha\ast F(u_a)\right)f(u_a)$, which proves
    \begin{equation*}
	\left(I_\alpha\ast F(u_n)\right)f(u_n)\rightharpoonup\left(I_\alpha\ast F(u_a)\right)f(u_a)\qquad\mbox{in}\ \,L^1(\R^2)\,.
    \end{equation*}
    This, together with $u_n\rightharpoonup u_a$ in $H^1(\R^2)$ and $\lambda_n\to\lambda_a$ yields \eqref{20240531-e6}, and in turns the Poho\v zaev identity \eqref{20240531-e8} for $(u_a,\lambda_a)$. Hence, by \eqref{20240531-e7} and \eqref{20240531-e8} it follows that $u_a\in\Sra$ and $u_n\to u_a$ in $L^2(\R^2)$, and moreover $\cP(u_a)=0$, namely
    \begin{align}\label{20240531-e9}
        \|\nabla u_a\|_2^2+\frac{2+\alpha}2\int_{\R^2}(I_\alpha\ast F(u_a))F(u_a)\dd x-\intd(I_\alpha\ast F(u_a))f(u_a)u_a\dd x=0\,.
    \end{align}
    By \eqref{20240531-e9} and ($f_3$),
    \begin{align}\label{20240531-e10}
        \frac12\|\nabla u_a\|_2^2-\frac12\intd(I_\alpha\ast F(u_a))F(u_a)\dd x&=\frac12\intd(I_\alpha\ast F(u_a))\left[f(u_a)u_a-\left(2+\frac\alpha2\right)F(u_a)\right]\dd x\notag\\
        &\geq\frac{\mu-(2+\frac\alpha 2)}{2}\intd(I_\alpha\ast F(u_a))F(u_a)\dd x\geq0\,.
    \end{align}
    Furthermore, by \eqref{20240531-e4} and the Brezis-Lieb Lemma, since $\{u_n\}$ is bounded in $H^1(\R^2)$, we have
    \begin{align}\label{20240531-e5}
        c_{mp}=\frac12\left(\lim_{n\to+\infty}\|\nabla(u_n-u_a)\|_2^2+\|\nabla u_a\|_2^2\right)-\frac12\int_{\R^2}(I_\alpha\ast F(u_a))F(u_a)\dd x\,.
    \end{align}
    Combining \eqref{20240531-e5} and \eqref{20240531-e10} we obtain
    \begin{align*}
        c_{mp}\geq\frac12\lim_{n\to+\infty}\|\nabla(u_n-u_a)\|_2^2\,,
    \end{align*}
    and so, by Corollary \ref{MP_level_bound}
    \begin{align*}
        \lim_{n\to+\infty}\|\nabla(u_n-u_a)\|_2^2<\frac{(2+\alpha)\pi}{\gamma_0}\,.
    \end{align*}
    According to Lemma \ref{lemma 2.3} we deduce that
    \begin{align}\label{20240531-e11}
        \lim_{n\to+\infty}\int_{\R^2}(I_\alpha\ast F(u_n))f(u_n)u_n\dd x=\int_{\R^2}(I_\alpha\ast F(u_a))f(u_a)u_a\dd x\,.
    \end{align}
    Now, combining \eqref{20240531-e4}, \eqref{20240531-e9}, \eqref{20240531-e11}, and the fact of $\cP(u_n)\to 0$, we have that $\|\nabla u_n\|_2^2\to\|\nabla u_a\|_2^2.$ Hence, $u_n\to u_a\not\equiv0$ in $H_r^1(\R^2)$.
\end{proof}

\section{Proof of Theorems \ref{Thm} and \ref{Thm_gs}}\label{sec_proofs}
\begin{proof}[Proof of Theorem \ref{Thm}]
    By Lemma \ref{20240605L2}, there exists a $(PSP)_{c_{mp}}$-sequence for $J$ constrained to $\Sra$. It follows from Lemma \ref{MP_geometry} that $c_{mp}>0$. Furthermore, by Corollary \ref{MP_level_bound} we have $c_{mp}<\frac{(2+\alpha)\pi}{2\gamma_0}$. Then the compactness Lemma \ref{lemma 3.6} applies and we can find a mountain pass type critical point $u_a\in H_r^1(\R^2)$ of $\left.J\right|_{\Sra}$, which weakly solves \eqref{eq} for some $\lambda_a>0$. Moreover, testing \eqref{eq} by $u_a^-:=-\min\{u_a,0\}\in H^1(\R^2)$, one can see that
    \begin{equation*}
        \begin{split}
            0&=J'(u_a)[u_a^-]+\lambda_aN'(u_a)[u_a^-]\\
            &=\intd\nabla u_a\nabla u_a^-\dd x-\intd\left(I_\alpha\ast F(u_a)\right)f(u_a)u_a^-\dd x+\lambda_a\intd u_a\,u_a^-\dd x\\
            &=\intd|\nabla u_a^-|^2\dd x+\lambda_a\intd|u_a^-|^2\dd x\,,
        \end{split}
    \end{equation*}
    by ($f_0$), hence $u_a^-\equiv0$; thus, $u_a\geq 0$ in $\R^2$. Finally, by the strong maximum principle for semilinear equations (see e.g. \cite[Theorem 11.1]{PucciSerrin}), we conclude that $u_a$ is positive in $\R^2$.

    We are left to prove the a priori estimate on the frequency $\lambda_a$. Exploiting the identities $J\big|_{\Sra}'(u_a)[u_a]-\cP(u_a)=0$, $\|u_a\|_2^2=a^2$, and $J(u_a)=c_{mp}$, one infers
    \begin{equation*}
        \begin{split}
            \lambda_a a^2=\frac{2+\alpha}2\intd\left(I_\alpha\ast F(u_a)\right)F(u_a)\dd x=\frac{2+\alpha}2\left(\|\nabla u_a\|_2^2-2c_{mp}\right).
        \end{split}
    \end{equation*}
    This, together with the uniform upper bound on $\|\nabla u_n\|_2$ proved in \eqref{bound_un}, the strong convergence $u_n\to u_a$ in $H^1(\R^2)$ from Lemma \ref{lemma 3.6}, and the bound for the mountain pass level found in Corollary \ref{MP_level_bound}, yields
    \begin{equation*}
        \begin{split}
            \lambda_a&=\frac{2+\alpha}{2a^2}\left(\|\nabla u_a\|_2^2-2c_{mp}\right)\leq\frac{2+\alpha}{2a^2}\left(\frac{2\mu-(2+\alpha)}{\mu-\left(2+\frac\alpha2\right)}-2\right)c_{mp}\\
            &=\frac{2+\alpha}{2a^2}\cdot\frac{2c_{mp}}{\mu-\left(2+\frac\alpha2\right)}\leq\frac{(2+\alpha)^2\pi}{2a^2\gamma_0\left(\mu-\left(2+\frac\alpha2\right)\right)}\,,
        \end{split}
    \end{equation*}
    which is the desired bound.
\end{proof}

We recall that the Poho\v zaev manifold $\cP(a)$ is defined in \eqref{Pohozaev_manifold}, and we denote by $\cP_r(a):=\{u\in\Sra\,|\,\cP(u)=0\}$ its radial counterpart. Moreover, we define the following least-energy levels
$$m(a):=\inf_{v\in\cP(a)}J(v)\,,\qquad m_r(a):=\inf_{v\in\cP_r(a)}J(v)\,.$$

\begin{Lem}\label{DYlemma4.6}
    Under ($f_0$)–($f_3$) and ($f_6$), if $u\in\cS_a$, then $s\mapsto\tJ(u,s)$ reaches its unique maximum at a point $s(u)\in\R$ such that $H(u,s(u))\in\cP(a)$.
\end{Lem}
\begin{proof}
    Here, we retrace the argument in \cite[Lemma 4.6]{DY}. By Lemmas \ref{H_J} and \ref{MP_geometry} we know that for an arbitrary $u\in\Sra$ there exists (at least) $s_0\in\R$ such that $\tfrac{\dd\tJ(u,s)}{\dd s}\big|_{s=s_0}=0$, hence $\cH(u,s_0)\in\cP(a)$. Computing the derivative, we have
    \begin{equation*}
        \frac{\dd\tJ}{\dd s}(u,s)=\e^{2s}\left(\intd|\nabla u|^2\dd x-\psi_u(s)\right),
    \end{equation*}
    where
    $$\psi_u(s):=\intd\left(I_\alpha\ast\frac{F\left(\e^su\right)}{(\e^s)^{2+\frac\alpha2}}\right)\frac{\tF\left(\e^su\right)}{(\e^s)^{2+\frac\alpha2}}\dd x\,.$$
    On the one hand, ($f_6$) implies that $\sigma\mapsto\frac{\tF\left(\sigma t\right)}{\sigma^{2+\frac\alpha2}}$ is nondecreasing for any $t\in\R$; on the other hand,
    $$\frac\dd{\dd\sigma}\left(\frac{F(\sigma t)}{\sigma^{2+\frac\alpha2}}\right)=\frac{f(\sigma t)\sigma t-\left(2+\frac\alpha2\right)F(\sigma t)}{\sigma^{3+\frac\alpha2}}>0\,.$$
    Hence the function $\psi_u$ is strictly increasing in $\R$. This implies that the stationary point is unique and is a maximum.
\end{proof}

\begin{Lem}\label{c2=c2rad}
    Under ($f_0$)–($f_3$) and ($f_6$), one has $m(a)=m_r(a)\,$.
\end{Lem}
\begin{proof}
    It is clear that $m(a)\leq m_r(a)$. We aim at showing the reverse inequality by means of a rearrangement argument. Note indeed that $u\in\cS_a$ implies $u^*\in\Sra$, where $u^*$ is its Schwarz rearranged function, since $\|u\|_2=\|u^*\|_2$. Let $u\in\cS_a$, by Lemma \ref{DYlemma4.6} there exist $s_1,s_2\in\R$ such that $\displaystyle{\max_{s\in\R}}\,\tJ(u,s)=\tJ(u,s_1)$ and $\displaystyle{\max_{s\in\R}}\,\tJ(u^*,s)=\tJ(u^*,s_2)$. Hence, by means of P\'olya-Szeg\H o inequality and the Riesz rearrangement inequality (see \cite[Theorem 3.7]{LiebLoss}),
    \begin{equation*}
        \begin{split}
            J(\cH(u,s_1))&\geq J(\cH(u,s_2))=\frac{\e^{2s_2}}2\intd|\nabla u|^2\dd x-\frac1{2\e^{(2+\alpha)s_2}}\intd\left(I_\alpha\ast F(\e^{s_2}u)\right)F(\e^{s_2}u)\dd x\\
            &\geq\frac{\e^{2s_2}}2\intd|\nabla u^*|^2\dd x-\frac1{2\e^{(2+\alpha)s_2}}\intd\left(I_\alpha\ast F(\e^{s_2}u^*)\right)F(\e^{s_2}u^*)\dd x=J(\cH(u^*,s_2))\,.
        \end{split}
    \end{equation*}
    Since $\cH(u^*,s_2)\in\cP_r(a)$, we deduce
    $$J(\cH(u,s_1))\geq J(\cH(u^*,s_2))\geq\inf_{v\in\cP_r(a)}J(v)=m_r(a)\,.$$
    The proof is concluded once we recall that $\cH(u,s_1)\in\cS_a$.
\end{proof}

\begin{proof}[Proof of Theorem \ref{Thm_gs}]
    In view of Lemma \ref{c2=c2rad}, all we need to prove is that $c_{mp}=m_r(a)$. Since $u_a\in\cP_r(a)$ by Theorem \ref{Thm}, being a normalised radial solution of \eqref{eq}, it is then clear that $c_{mp}\geq m_r(a)$. On the other hand, we claim that
    \begin{equation}\label{levels_equalities}
        c_{mp}\leq c_r(a):=\inf_{v\in\Sra}\max_{s\in\R}J(\cH(v,s))=m_r(a)\,.
    \end{equation}
    Following \cite{Dou2023}, take $v\in\Sra$, then by Lemma \ref{H_J} there exists $s_0>0$ such that $J(\cH(v,-s_0))\in A$ and $J(\cH(v,s_0))<0$. Let $s\mapsto\overline\cH(v,s):=\cH(v,(2s-1)s_0)$, then $\overline\cH(v,\cdot)\in\Gamma$, therefore
    $$\max_{s\in\R}J(\cH(v,s))=\max_{s\in[-s_0,s_0]}J(\cH(v,s))=\max_{s\in[0,1]}J(\overline\cH(v,s))\geq c_{mp}\,,$$
    from which one easily gets $c_{mp}\leq c_r(a)$. Finally, being trivial the fact that $c_r(a)\geq m_r(a)$, we need to show the reverse inequality. Supposing the strict inequality, we would find $v\in\cP_r(a)$ such that $J(v)<c_r(a)$. However,
    $$J(v)=J(\cH(v,0))=\max_{s\in\R}J(\cH(v,s))\geq c_r(a)\,.$$
    This justifies \eqref{levels_equalities} and the proof is concluded.
\end{proof}

\vskip0.4truecm

\paragraph{\textbf{Acknowledgements}:} The authors would like to express their gratitude to the referee for carefully reading the manuscript and providing numerous valuable comments that have significantly improved the paper. Ling Huang is supported by the Scientific Research Innovation Project of Graduate School of South China Normal University. Giulio Romani is a member of \textit{Gruppo Nazionale per l'Analisi Matematica, la Probabilità e le loro Applicazioni} (GNAMPA) of the \textit{Istituto Nazionale di Alta Matematica} (INdAM) and is partially supported by INdAM-GNAMPA Project 2024 titled \textit{New perspectives on Choquard equation through PDEs with local sources} (CUP E53C2200l93000l) and INdAM-GNAMPA Project 2025 titled \textit{Critical and limiting phenomena in nonlinear elliptic systems} (CUP E5324001950001).
\vskip0.4truecm
\paragraph{\textbf{Conflict of interest}:} There are no interests of a financial or personal nature.

\vskip0.4truecm

\end{document}